\RequirePackage{etoolbox}
\csdef{input@path}{{style/}{graphics/}}
\documentclass[MSNbibl,number,citesort,seceqn,rotating,dvips]{arxbj}
\usepackage{url,breakurl}
\usepackage{dcolumn}
\usepackage{upgreek}
\usepackage{graphicx}

\aid{0}
\volume{21}
\issue{1}
\pubyear{2015}
\firstpage{242}
\lastpage{275}
\doi{10.3150/13-BEJ566} 
\docsubty{FLA}

\makeatletter
\newcommand{\rrvert}{\vert}
\newcommand{\llvert}{\vert}
\newproclaim{defi}{Definition}[section]
\newtheorem{theo}{Theorem}[section]
\newtheorem{coll}{Corollary}[section]
\newremark{exmp}{Example}[section]
\newremark{Remarks}{Remarks}

\newcolumntype{d}[1]{D{.}{.}{#1}}
\newcommand{\boldintav}{\mathbf{v}}
\renewcommand{\epsilon}{\varepsilon}
\newcommand{\A}{\mathbf{A}}
\newcommand{\B}{\mathbf{B}}
\newcommand{\C}{\mathbf{C}}
\newcommand{\D}{\mathbf{D}}
\renewcommand{\H}{\mathbf{H}}
\newcommand{\I}{\mathbf{I}}
\renewcommand{\O}{\mathbf{O}}
\renewcommand{\P}{\mathbf{P}}
\newcommand{\bfitV}{\mathbf{V}}
\newcommand{\bfitX}{\mathbf{X}}
\newcommand{\sbfitX}{\mathbf{X}}
\newcommand{\bfitb}{\mathbf{b}}
\newcommand{\bfite}{\mathbf{e}}
\newcommand{\bfitu}{\mathbf{u}}
\newcommand{\bfitw}{\mathbf{w}}
\newcommand{\bfitz}{\mathbf{z}}
\newcommand{\alp}{\boldsymbol{\alpha}}
\newcommand{\be}{\boldsymbol{\beta}}
\newcommand{\bfbeta}{\boldsymbol{\beta}}
\newcommand{\bfGam}{\boldsymbol{\Gamma}}
\newcommand{\bfPi}{\boldsymbol{\Pi}}
\makeatother

\begin{document}
\begin{frontmatter}

\title{Variable selection and estimation for semi-parametric
multiple-index models}
\runtitle{Estimation for multi-index models}

\begin{aug}
\author[1]{\inits{T.}\fnms{Tao}~\snm{Wang}\thanksref{1}},
\author[2]{\inits{P.}\fnms{Peirong}~\snm{Xu}\thanksref{2}} 
\and
\author[1]{\inits{L.}\fnms{Lixing}~\snm{Zhu}\corref{}\thanksref{1}\ead[label=e3]{lzhu@hkbu.edu.hk}}
\address[1]{Department of Mathematics, Hong Kong Baptist University, Hong Kong, China.\\ \printead{e3}}
\address[2]{Department of Mathematics, Southeast University, Nanjing, China}
\end{aug}

\received{\smonth{9} \syear{2012}}
\revised{\smonth{5} \syear{2013}}

%
\begin{abstract}
In this paper, we propose a novel method to select significant
variables and estimate the corresponding coefficients in multiple-index
models with a group structure. All existing approaches for single-index
models cannot be extended directly to handle this issue with several
indices. This method integrates a popularly used shrinkage penalty such
as LASSO with the group-wise minimum average variance estimation. It is
capable of simultaneous dimension reduction and variable selection,
while incorporating the group structure in predictors. Interestingly,
the proposed estimator with the LASSO penalty then behaves like an
estimator with an adaptive LASSO penalty. The estimator achieves
consistency of variable selection without sacrificing the root-$n$
consistency of basis estimation. Simulation studies and a real-data
example illustrate the effectiveness and efficiency of the new method.
\end{abstract}

%
\end{frontmatter}

\section{Introduction}\label{sec1ch5}

Suppose that $Y \in\mathbb{R}$ is a univariate response and $\bfitX=
(X_1, \ldots, X_p)^{\top} \in\mathbb{R}^p$ is a vector of
predictors. A general goal of regression analysis is to characterize
the conditional distribution of $Y$ given $\bfitX$, or the conditional
mean $E(Y|\bfitX)$. The theory of sufficient dimension reduction (Li
\cite{Li91N1} and Cook and Weisberg \cite{b10}) provides a framework for reducing the
dimension of $\bfitX$ while preserving information on regression. Let
$\mathcal{S}$ denote a subspace of $\mathbb{R}^p$, and let $\P
_{\mathcal{S}}$ denote the orthogonal projection onto $\mathcal{S}$
with respect to the usual inner product. If $Y$ and $\bfitX$ are
independent conditioned on $\P_{\mathcal{S}}\bfitX$, then we say
that $\mathcal{S}$ is a dimension reduction subspace. The intersection
of all such subspaces, if itself
satisfies the conditional independence, is defined to be the central
subspace (Cook \cite{Coo98} and Yin, Li and Cook \cite{YinLiCoo08}). When only the mean
response $E(Y|\bfitX)$ is of interest, sufficient dimension reduction
can be defined in a similar fashion. Specifically, a subspace $\mathcal
{S}$ is said to be a mean dimension reduction subspace if $Y$ is
independent of $E(Y|\bfitX)$ given $\P_{\mathcal{S}}\bfitX$. If the
intersection of all mean dimension reduction subspaces is also a mean
dimension reduction subspace, it is called the central mean subspace
(Cook and Li \cite{CooLi02}). In either case, sufficient
dimension reduction
permits us to restrict attention to a number $d \leq p$ of new
predictors, expressed as linear combinations of the original ones:
$\bfbeta_1^{\top} \bfitX, \ldots, \bfbeta_d^{\top} \bfitX$,
where $\{\bfbeta_1, \ldots, \bfbeta_d\}$ is a basis of $\mathcal{S}$.

In the last two decades or so, a series of papers have considered
issues related to dimension reduction in regression. There have
primarily been two
categories of estimation methods in the literature: inverse
regression methods (Li \cite{Li91N1}, Cook and Weisberg \cite{b10} and Cook and Ni \cite
{CooNi05}) and direct regression methods (H\"{a}rdle and Stoker \cite{HarSto89},
Xia \textit{et al.} \cite{Xiaetal02} and Yin and Li \cite{YinLi11}). Inverse
regression methods, despite being
computationally simple and widely used, require strong assumptions on
predictors such as the linearity condition (Li \cite{Li91N1}), and often fail to
estimate the central subspace exhaustively (Cook
\cite{Coo98}). In contrast, the minimum average variance
estimation (MAVE)
method of Xia \textit{et al.} \cite{Xiaetal02} has proven effective in dimension reduction
and estimation of
complicated semi-parametric models. The
root-$n$ consistency is still achievable for the MAVE estimate.
Compared with other direct regression methods, the calculation for
MAVE is much easier, and many efficient algorithms are available.
Although MAVE was originally proposed for dimension reduction for the
conditional mean, the idea was recently generalized to target the
central subspace (Wang and Xia \cite{WanXia08} and
Yin and Li \cite{YinLi11}). In this article,
we are concerned mainly with predictors in the conditional mean.

Dimension reduction is a fundamental statistical problem in both theory
and practice. The aforementioned dimension reduction methods, however,
suffer from the difficulty of interpreting the results, because the new
extracted predictors usually involve all of the original ones. To
handle this problem, model-free variable selection, in the framework of
sufficient dimension reduction, has attracted considerable attention in
recent years. For example, Li, Cook and Nachtsheim \cite{LiCooNac05} introduced
test-based procedures, Bondell and Li \cite{BonLi09}
incorporated inverse
regression estimation with LASSO (Tibshirani \cite{Tib96}) to
obtain shrinkage
inverse regression estimation, and Chen, Zou and Cook \cite{CheZouCoo10} proposed a
unified method called coordinate-independent sparse estimation. See
also Zhu \textit{et al.} \cite{Zhuetal11}  and Wang, Xu and Zhu \cite{WanXuZhu12}. All these methods,
which are largely ``parametric'' in nature, are based on inverse
regression methods and thus suffer the drawbacks of strong design
assumptions and poor finite-sample performance (Wang, Xu and Zhu
\cite{WanXuZhu13}).

Exploring the idea of combining MAVE and LASSO, Wang and Yin
\cite{WanYin08}
proposed a sparse MAVE method and Zeng, He and Zhu \cite{ZenHeZhu12}
designed for single-index models a lasso-type approach called
sim-lasso. Because the sparse MAVE penalizes the index vectors
directly, it is not a principled method for variable selection and only
provides a sparse estimate for a basis matrix of the central mean
subspace column by column. The use of the $l_1$ penalty function in
Zeng, He and Zhu \cite{ZenHeZhu12} is novel in that it
penalizes the index vector
and the norm of the derivative of link function simultaneously.
However, the theoretical properties of sim-lasso, such as its
consistency and convergence rate, have not yet been studied due to the
interaction between the bandwidth and the penalty parameter. Further,
it is nontrivial, if not impossible, to extend sim-lasso to deal with
multiple-index models. Several papers have addressed the problem of
semi-parametric variable selection for single-index models, and
developed large sample properties. See, for instance, Liang \textit{et al.}
\cite{Liaetal10}, Peng and Huang \cite{PenHua11} and Wang, Xu and Zhu \cite{WanXuZhu13}. However,
condition (vi) in Liang \textit{et al.} \cite{Liaetal10} may not hold true and their
approach could not be extended to handle multiple-index models. The
penalized MAVE method in Wang, Xu and Zhu \cite{WanXuZhu13} was motivated by the
reasoning that predictor selection can be realized through selection of
nonvanishing rows of a basis matrix of the central mean subspace. A
bridge penalty function was employed to penalize the $l_1$ norms of the
rows of a basis matrix. Although the penalized MAVE performs well for
multiple-index models in the numerical studies, its theoretical
properties are established only for the special case of single-index
models. This is because, condition (C5) in Wang, Xu and Zhu
\cite{WanXuZhu13},
which is also assumed in Peng and Huang \cite{PenHua11},
is hard to check and
possibly invalid except for single-index models. To the best of our
knowledge, semi-parametric variable selection for multiple-index models
has thus far not been well studied.

In many engineering and scientific situations, however, predictors are
naturally grouped. For example, in biological applications assayed
genes or proteins can be grouped by biological pathways. Although
useful, existing dimension reduction methods are generic and treat all
predictors in $\bfitX$ indiscriminately. To take advantage of such
group knowledge, Li, Li and Zhu \cite{LiLiZhu10}
proposed a group-wise sufficient
dimension reduction method, called group-wise MAVE, which preserves
full regression information in the conditional mean of $Y$ given
$\bfitX$ while exploiting the group structure among predictors.
Generally, it is believed that incorporating group information into
dimension reduction can facilitate interpretation of results and
improve estimation accuracy as the number of unknown parameters has
been greatly reduced.

As a simple illustration, we use an example to show the necessity of
group-wise dimension reduction and variable selection. Consider a
response model $Y =
\be_1^{\top}\bfitX_1 + \sin(0.2\uppi\be_2^{\top}\bfitX_2) + 0.5
\epsilon$,
where $\bfitX_1 \in\mathbb{R}^{10}$, $\bfitX_2 \in\mathbb{R}^{10}$,
$\be_1 =
(1, -1, 0, \ldots, 0)^{\top}$, $\be_2 = (1, 1, 0, \ldots, 0)^{\top}$,
and all predictors and $\epsilon$ are independent standard normal
variables. Write $\bfitX= (\bfitX_1^{\top}, \bfitX_2^{\top})^{\top
}$. Then the
central mean subspace for $E(Y|\bfitX)$ is spanned by $(\be_1^{\top},
\mathbf{0}_{10}^{\top})^{\top}$ and $(\mathbf{0}_{10}^{\top},
\be_2^{\top})^{\top}$, where $\mathbf{0}_{10}$ is a $10
\times1$
vector of zeros. We then should rule out zeros and identify
$\be_1$ and $\be_2$ or their linear combinations. A~single but representative
simulated data set with 150 observations was obtained, and the MAVE
direction estimates were
\begin{eqnarray*}
&&(-0.624, 0.647, 0.006, 0.057, 0.034, -0.010, -0.013, 0.033, 0.023, -0.022,
\\
&&\quad{}-0.275, -0.316, -0.002, -0.062, -0.057, 0.010, 0.005, -0.034, -0.028,
-0.017)^{\top}
\end{eqnarray*}
and
\begin{eqnarray*}
&&(0.141, 0.379, -0.394, 0.005, 0.030, -0.313, -0.313, 0.146, 0.201, -0.341,
\\
&&\quad 0.106, -0.022, -0.286, 0.111, 0.096, -0.047, -0.303, -0.073, -0.203,
-0.226)^{\top}.
\end{eqnarray*}
MAVE treats all predictors in $\bfitX$ indiscriminately. While the
first direction estimate seems reasonable, the second one
is very poor, and thus the overall estimation accuracy must be poor.
Given the prior information that
$\bfitX_1$ and $\bfitX_2$ are two predictor groups, we apply group-wise
MAVE, and the resulting estimates of $\be_1$ and $\be_2$ are
respectively, given by
\[
(0.698, -0.715, -0.035, -0.031, -0.010, -0.003, -0.001, 0.010, -0.010,
0.015)^{\top}
\]
and
\[
(0.717, 0.666, -0.035, 0.155, 0.002, 0.032, -0.008, 0.111, -0.024,
-0.058)^{\top}.
\]
%
A substantial gain in accuracy has been
achieved by incorporating the predictor group information.
Nevertheless, in each group all the
predictors are included in the extracted linear combination, although
some coefficients are small, obscuring the fact that only the first two
predictors are contributing factors. It is obvious that group-wise MAVE
cannot be the base for both dimension reduction and variable
selection. Therefore, a selection operator also plays an important
role, and we will see that a shrinkage penalty will be useful for us to
use group-wise MAVE to exclude irrelevant predictors from the model.

Two main features of this paper are listed below.
\begin{enumerate}[2.]
\item We consider the problem of semi-parametric variable selection for
multiple-index regression models. Although multiple-index models are
popular in the statistics and econometrics literature, little work has
been done on variable selection.
We propose a shrinkage
MAVE estimator by introducing a shrinkage factor for each row of an
estimated basis matrix of the central mean subspace. For multiple-index
models the proposed estimator is proved to be consistent in variable
selection while retaining the root-$n$ consistency. However, although
the estimation problem can be reformulated as a LASSO problem in
spirit, the LASSO problem under study has an asymptotically singular
design matrix (Knight and Fu \cite{KniFu00}). This is
because the MAVE procedure
is a combination of nonparametric function estimation and direction
estimation. This makes the theoretical investigation more complicated.
To deal with this issue for single-index models, condition (C5) is
assumed in Wang, Xu and Zhu \cite{WanXuZhu13},
otherwise, the large sample
properties are difficult to derive. For multiple-index models, the
standard approach of LASSO with nonsingular designs fails to show the
large sample properties. Therefore, in this paper, the results of
mixed-rates asymptotics (Radchenko \cite{Rad08}) are
adopted to derive the
asymptotic behavior even the design matrix is asymptotically singular.
This is a new skill about proving the asymptotics of the LASSO
estimation for semi-parametric models. The interaction between the
bandwidth and the penalty parameter now is explicitly shown in Theorem~\ref{theo1ch5}.
%
\item We propose a general knowledge-based method that accounts for
prior group information. As we have explained before, the group
structure leads to a reduction in the total number of parameters.
Consequently, our method, which is motivated by and derives from
dimension reduction, doubly alleviates the ``curse of dimensionality''.
As a by product, such a structure also makes the computation more efficient.
\end{enumerate}
%

The paper is organized as follows. In Section~\ref{sec21ch5}, we
review the group-wise minimum average variance estimation. In
Section~\ref{sec22ch5}, we combine group-wise MAVE with the LASSO
penalty, as
an example, to propose a shrinkage group-wise MAVE estimator. This
method does not require any restrictive design assumptions, and is
capable of simultaneous dimension reduction and variable selection. The
asymptotic properties of the new estimator are established in
Section~\ref{sec23ch5}. We also use a criterion, which has the same
form as the Bayesian information criterion (BIC; Schwarz
\cite{Sch78}), to
select the optimal tuning parameter. Moreover, we establish the
consistency of the resulting BIC-type selector. Numerical studies are
presented in Section~\ref{sec3ch5}. As many shrinkage penalties can
also be applied, we then include the simulation results with two other
penalties as well. All technical proofs are relegated
to the \hyperref[app]{Appendix}.


\section{Methodology}\label{sec2ch5}

We begin with some basic notations and terminology. For a positive
integer $m$, $\I_m$ stands for the $m \times m$ identity matrix. For
an $m_1 \times m_2$ matrix $\A$, $\operatorname{span}(\A)$
represents the column space of $\A$ and $\P_{\A}$ represents
the orthogonal projection onto $\operatorname{span}(\A)$. For a
subspace $\mathcal{S}$ of $\mathbb{R}^m$, if $\A$ is a matrix of
full column rank and $\operatorname{span}(\A) = \mathcal{S}$, then
we call $\A$ a basis matrix of $\mathcal{S}$. Moreover, $\P
_{\mathcal{S}}$ represents the projection onto $\mathcal{S}$, that
is, $\P_{\mathcal{S}} = \P_{\A}$, where $\A$ is any basis
matrix of $\mathcal{S}$. For an $m$-dimensional vector $\bfitw= (w_1,
\ldots, w_m)^{\top}$, $\operatorname{diag}(\bfitw)$ denotes a
diagonal matrix whose diagonal entries starting in the upper left
corner are $w_1, \ldots, w_m$. We use $\A_1 \oplus\cdots\oplus{\A
_g}$, or simply $\bigoplus_{l = 1}^g \A_l$, to denote a block diagonal
matrix with matrices $\A_1, \ldots, \A_g$ on the diagonal.

\subsection{A short review}\label{sec21ch5}

In this subsection, we review group-wise dimension reduction for the
regression mean function and the group-wise minimum average variance
estimation. We refer the reader to Li, Li and Zhu \cite{LiLiZhu10} for more details.

Let $\mathcal{S}_1, \ldots, \mathcal{S}_g$ be subspaces of
$\mathbb{R}^p$ that form an orthogonal decomposition of
$\mathbb{R}^p$, that is, $\mathbb{R}^p = \mathcal{S}_1\oplus\cdots
\oplus\mathcal{S}_g$, where $\oplus$ denotes the direct sum
operator. If there are subspaces $\mathcal{T}_l \subseteq
\mathcal{S}_l$ for $l = 1, \ldots, g$ such that $E(Y|\bfitX) =
E(Y|\P_{\mathcal{T}_1}, \ldots, \P_{\mathcal{T}_g})$, then we say
that $\mathcal{T}_1 \oplus\cdots\oplus\mathcal{T}_g$ is a
group-wise mean dimension reduction subspace with respect to
$\{\mathcal{S}_1, \ldots, \mathcal{S}_g\}$. Under very mild
conditions (Yin, Li and Cook \cite{YinLiCoo08}), the
intersection of all group-wise
mean dimension reduction subspaces, with respect to a given
orthogonal decomposition $\{\mathcal{S}_1, \ldots, \mathcal{S}_g\}$,
exists uniquely. We call this subspace the group-wise central mean
subspace and denote it as
$\mathcal{S}_{E(Y|\sbfitX)}(\mathcal{S}_1, \ldots,
\mathcal{S}_g)$. By definition,
\[
\mathcal{S}_{E(Y|\sbfitX)}(\mathcal{S}_1, \ldots,
\mathcal{S}_g) = \mathcal{T}_1^* \oplus\cdots\oplus
\mathcal{T}_g^*
\]
for some subspaces $\mathcal{T}_1^* \subseteq\mathcal{S}_1, \ldots,
\mathcal{T}_g^* \subseteq\mathcal{S}_g$. Let $p_l, d_l$ and $d$
denote the dimensions of $\mathcal{S}_l$, $\mathcal{T}_l^*$ and
$\mathcal{S}_{E(Y|\sbfitX)}(\mathcal{S}_1, \ldots, \mathcal
{S}_g)$, respectively. Then we have $p = p_1 + \cdots+ p_g$ and $d =
d_1 + \cdots+ d_g$.

Let $\bfGam_l \in\mathbb{R}^{p \times p_l}$ be a basis matrix of
$\mathcal{S}_l$, and let $\bfitV_l = \bfGam_l^{\top} \bfitX\in
\mathbb{R}^{p_l}$. We note that components of $\bfitV_l$ correspond
to predictors in group $l$, and all the group information contained in
$\bfGam_l$'s is available as prior knowledge. By construction, there
are matrices $\B_l^* \in\mathbb{R}^{p_l \times d_l}$ for $l = 1,
\ldots, g$ such that $\operatorname{span}(\bfGam_l\B_l^*) =
\mathcal{T}_l^*$. Write $\B^* = \bigoplus_{l = 1}^g\B_l^*$. We are
interested in estimating $\B^*$ or its column space $\operatorname
{span}(\B^*)$.

Li, Li and Zhu \cite{LiLiZhu10} proposed the group-wise
MAVE estimator such that
the matrix $\B^* = \bigoplus_{l = 1}^g \B_l^*$ is the minimizer of
\[
E\bigl\{Y - E\bigl(Y|\B_1^{\top}\bfitV_1,
\ldots, \B_g^{\top}\bfitV_g\bigr)\bigr
\}^2
\]
with respect to ${\B}_1 \in\mathbb{R}^{p_1 \times d_1}, \ldots, {\B
}_g \in\mathbb{R}^{p_g \times d_g}$, subject to $\B_l^{\top}\B_l =
\I_{d_l}$ for $l = 1, \ldots, g$.

Let $\bfitV= (\bfitV_1^{\top}, \ldots, \bfitV_g^{\top})^{\top}
\in\mathbb{R}^p$ and $\B= \bigoplus_{l = 1}^g \B_l$. Then we have
\[
E\bigl\{Y - E\bigl(Y|\B_1^{\top}\bfitV_1,
\ldots, \B_g^{\top}\bfitV_g\bigr)\bigr
\}^2 = E\bigl\{\sigma_{\B}^2\bigl(
\B^{\top}\bfitV\bigr)\bigr\},
\]
where
$\sigma_{\B}^2(\B^{\top}\bfitV) = E[\{Y - E(Y|\B_1^{\top
}\bfitV_1, \ldots, \B_g^{\top}\bfitV_g)\}^2|\B_1^{\top}\bfitV
_1, \ldots, \B_g^{\top}\bfitV_g]$ is the conditional variance of
$Y$ given $\B_1^{\top}\bfitV_1, \ldots, \B_g^{\top}\bfitV_g$.

Suppose that $\{(y^i, \boldintav ^i), i = 1, \ldots,
n\}$ is a random sample from $(Y, \bfitV)$. Extending the MAVE idea,
we can use local linear smoothing to estimate
$\sigma_{\B}^2(\B^{\top}\bfitV)$. Specifically, for
any given $\boldintav ^0 \in\mathbb{R}^p$, we have the following
approximation
\begin{eqnarray*}
\sigma_{\B}^2\bigl(\B^{\top}\boldintav
^0\bigr) &\approx& \sum_{j = 1}^n
\bigl\{y^j - E\bigl(Y|\B^{\top} \bfitV= \B^{\top}
\boldintav ^j\bigr)\bigr\}^2 w_{j}^0
\\
&\approx& \sum_{j = 1}^n \Biggl
\{y^j - a^0 - \sum_{l = 1}^g
\bfitb _{l}^{0\top} \B_l^{\top}\bigl(
\boldintav _l^j - \boldintav _{l}^0
\bigr) \Biggr\}^2 w_{j}^0,
\end{eqnarray*}
%
where $w_{j}^0$'s are kernel weights such that $\sum_{j
= 1}^{n}w_{j}^0 = 1$, and $a^0 + \sum_{l = 1}^g\bfitb_{l}^{0\top} \B
_l^{\top}(\boldintav _l^j - \boldintav _{l}^0)$ is the local linear expansion of $E(Y|\B^{\top} \bfitV
= \B^{\top} \boldintav ^j)$ at $\boldintav ^0$.

Consequently, we can recover the group-wise central mean subspace by
minimizing the objective function
%
%
\begin{equation}
\label{2.1.1ch5} \sum_{i = 1}^n\sum
_{j = 1}^n \Biggl\{y^j -
{{a}^i} - \sum_{l =
1}^g
\bfitb_{l}^{i\top} \B_l^{\top}\bigl(
\boldintav _l^j - \boldintav _{l}^i
\bigr) \Biggr\}^2 w_{j}^i
\end{equation}
with respect to ${a}^i \in\mathbb{R}$, $\bfitb_1^i \in
\mathbb{R}^{d_1}, \ldots, \bfitb_g^i \in\mathbb{R}^{d_g}$, $i = 1,
\ldots, n$, and $\B_l \in\mathbb{R}^{p_l \times d_l}$ with $\B
_l^{\top}\B_l = \I_{d_l}$ for $l = 1, \ldots, g$. To allow the
estimation to be adaptive to the regression structure, we follow the
idea of refined MAVE (Xia \textit{et al.} \cite{Xiaetal02} and Li, Li and Zhu
\cite{LiLiZhu10}) and adopt
the weights
\[
w_{j}^i = \frac{K_h\{\B^{\top}(\boldintav ^j -
\boldintav ^i)\}}{\sum_{j = 1}^n K_h\{\B^{\top
}(\boldintav ^j - \boldintav ^i)\}},
\]
where $K_h(\cdot)$ is a $d$-dimensional kernel with bandwidth $h$, and
$\B$ is taken to be the current or latest estimate.

The minimization problem in (\ref{2.1.1ch5}) can be solved by fixing
$({a}^i, \bfitb_1^i, \ldots, \bfitb_g^i)$, $i = 1, \ldots, n$, and
fixing $\{\B_l\}_{l = 1}^g$ alternatively. Thus, the calculation can
be decomposed into two optimization problems both of which have simple
analytic solutions. The details of the group-wise MAVE algorithm can be
found in Section~3.2 of Li, Li and Zhu \cite{LiLiZhu10}. Let $\tilde{\B} =
\bigoplus_{l = 1}^g \tilde{\B}_l$ denote the group-wise minimum average
variance estimator.

\subsection{Shrinkage group-wise minimum average variance
estimation}\label{sec22ch5}

The group-wise MAVE method captures the full regression information in
$E(Y|\bfitX)$ while preserving the group structure in $\bfitX$.
Specifically, it can provide a consistent estimator of $\B^*\D^0$ for
some $d \times d$ nonsingular matrix $\D^0 = \bigoplus_{l = 1}^g \D
_{l}^0$, where $\D_l^0 \in\mathbb{R}^{d_l\times d_l}$ for $l = 1,
\ldots, g$. However, the elements of $\tilde{\B}_l$'s are usually
nonzero. Consequently, the extracted predictor vector $\tilde{\B
}_l^{\top}\bfitV_l$ corresponding to group $l$ consists of linear
combinations of all the predictors in that group. When there are a
large number of predictors, one would expect that only a subset of
predictors are relevant to the response variable. Write $\bfitV= (V_1,
\ldots, V_p)^{\top}$. According to Proposition~1 of Cook \cite{Coo04},
$V_s$ is irrelevant if and only if the $s$th row of $\B^*$ is a zero
vector. Further, it is easy to see that for any $d \times d$
nonsingular matrix $\D$, when a row of $\B^*$ is zero, the
corresponding row of $\B^*\D$ is also zero, and vice versa. These
observations motivate us to employ the state-of-the-art methods for
simultaneous shrinkage estimation and variable selection, such as
LASSO, to design a sparse version of the group-wise MAVE procedure
which shrinkages some rows of $\tilde{\B}$ to be exactly zero vectors.

Define
\[
\tilde{w}_{j}^i = \frac{K_h\{\tilde{\mathbf{B}}^{\top}({\boldintav ^j} - {\boldintav }^i)\}}{\sum_{j = 1}^n K_h\{\tilde{\mathbf{B}}^{\top}({\boldintav }^j - {\boldintav }^i)\}}, \qquad i, j = 1,
\ldots, n.
\]
For each $i = 1, \ldots, n$, let $({\tilde{a}^i}, \tilde{\bfitb
}_1^i, \ldots, \tilde{\bfitb}_g^i)$ be the minimizer of
%
%
\begin{equation}
\sum_{j = 1}^n \Biggl\{y^j -
{{a}^i} - \sum_{l = 1}^g{\bfitb
}_{l}^{i\top} \tilde{\mathbf {B}}_l^{\top}
\bigl({\boldintav }_l^j - {\boldintav
}_{l}^i\bigr) \Biggr\}^2 \tilde{w}_{j}^i.
\end{equation}
In the sequel, we shall use an updated version of the group-wise
minimum average variance estimator, $\tilde{\tilde{\B}} = \bigoplus_{l
= 1}^g \tilde{\tilde{\B}}_l$, which is the minimizer of
%
%
\begin{equation}
\sum_{i = 1}^n\sum
_{j = 1}^n \Biggl\{y^j - {
\tilde{a}^i} - \sum_{l
= 1}^g
\tilde{\bfitb}_{l}^{i\top} \B_l^{\top}
\bigl(\boldintav _l^j - \boldintav _{l}^i
\bigr) \Biggr\}^2 \tilde{w}_{j}^i.
\end{equation}

%
\begin{defi}\label{defi1ch5}
A shrinkage group-wise minimum average variance estimator is defined as
\[
\hat{\mathbf {B}} = \bigoplus_{l = 1}^{g}
\operatorname{diag}(\hat{\alp }_l)\tilde{\tilde{\mathbf{B}}}_l,
\]
where the shrinkage index vectors $\hat{\alp}_l = (\hat{\alpha
}_{l1}, \ldots, \hat{\alpha}_{lp_l})^{\top} \in\mathbb{R}^{p_l}$
for $l = 1, \ldots, g$ are determined by minimizing
%
%
\begin{equation}
\label{2.2.1ch5} \sum_{i = 1}^{n}\sum
_{j = 1}^{n} \Biggl[y^j - {
\tilde{a}^i} - \sum_{l = 1}^{g}
\tilde{{\bfitb}}_l^{i\top} \bigl\{\operatorname{diag}(
\alp_l)\tilde {\tilde{\mathbf{B}}}_l\bigr
\}^{\top}\bigl({\boldintav }_l^j - {\boldintav
}_l^i\bigr) \Biggr]^2 \tilde{w}_{j}^i
\end{equation}
with respect to $\alp_l = ({\alpha}_{l1}, \ldots, {\alpha
}_{lp_l})^{\top} \in\mathbb{R}^{p_l}$, $l = 1, \ldots, g$, subject to
$\sum_{l = 1}^g\sum_{s = 1}^{p_l}|\alpha_{ls}| \leq\tau_n$ for
some $\tau_n \geq0$.
\end{defi}

To solve the above optimization problem, we note that (\ref{2.2.1ch5})
can be re-expressed as
\[
\sum_{i = 1}^{n}\sum
_{j = 1}^{n} \Biggl\{y^j - {
\tilde{a}^i} - \sum_{l = 1}^{g} {
\tilde{\bfitb}}_l^{i\top} \tilde{\tilde{\B
}}_l {\vphantom{\tilde{B}}\!\!}^{\top} \operatorname{diag}\bigl(\boldintav
_l^j - \boldintav _l^i\bigr)
\alp_l \Biggr\}^2 \tilde{w}_{j}^i.
\]
Equivalently, the shrinkage index vectors minimize
%
%
\begin{equation}
\label{2.2.2ch5} \sum_{i = 1}^{n}\sum
_{j = 1}^{n} \Biggl\{y^j - {
\tilde{a}^i} - \sum_{l = 1}^{g}{
\tilde{\bfitb}}_l^{i\top} \tilde{\tilde{\B}}_l{\vphantom{\tilde{B}}\!\!}^{\top
}
\operatorname{diag}\bigl(\boldintav _l^j - \boldintav
_l^i\bigr)\alp_l \Biggr\}^2
\tilde{w}_{j}^i + \lambda_n \sum
_{l = 1}^g\sum_{s = 1}^{p_l}|
\alpha_{ls}|
\end{equation}
for some tuning parameter $\lambda_n \geq0$. As a result,
commonly-used LASSO algorithms, such as those of Efron \textit{et al.} \cite{Efretal04}
and Friedman, Hastie and Tibshirani \cite{FriHasTib10}, can be applied
to obtain
the shrinkage index vectors $\hat{\alp}_l$ for $l = 1, \ldots, g$.

When $\tau_n \geq p$, the indices $\hat{\alpha}_{ls} = 1$ for all $l
= 1, \ldots, g$ and $s = 1, \ldots, p_l$, and so $\hat{\B}$ reduces
to the usual group-wise MAVE estimator $\tilde{\tilde{\B}}$. As
$\tau_n$ gradually decreases, some of the indices are shrunk to zero,
which means some rows of $\hat{\B}$ are zero; that is, the
corresponding predictors are irrelevant to the response variable given
the other predictors.

\subsection{Asymptotic theory}\label{sec23ch5}

We next study the large-sample properties of the proposed method. For
an $m_1 \times m_2$ matrix $\A$, we say that $\A$ is row-sparse if
some of its rows are zero. Let $\mathcal{I}(\A) \subseteq\{1, \ldots
, m\}$ denote the subset of indices corresponding to nonzero rows of
$\A$. Clearly, the notion of row-sparseness is
nonsingular-transformation independent, since for any $m_2 \times m_2$
nonsingular matrix $\O$, $\mathcal{I}(\A) = \mathcal{I}(\A\O)$.
Suppose that $\B^* = \bigoplus_{l = 1}^{g} \B_l^*$ is row-sparse.
Without loss of generality, we assume that for $l = 1, \ldots, g$ the
first $q_l$ rows of $\B_l^*$ are nonzero, that is, $\mathcal{I}(\B
_l^*) = \{1, \ldots, q_l\}$. The following theorem concerns the
asymptotic behavior of shrinkage group-wise MAVE.

%
%
\begin{theo}\label{theo1ch5}
Suppose that the regularity conditions \textup{(A1)--(A6)} given in the \hyperref[app]{Appendix}
hold. If $\lambda_n \rightarrow\infty$ and $\lambda_n
n^{-1/2}h^{-2} \rightarrow0$, then we have
\begin{longlist}[(2)]
\item[(1)] selection consistency: $P\{\mathcal{I}(\hat{\mathbf
{B}}_l) = \mathcal{I}(\mathbf {B}_l^*), l = 1, \ldots, g\} \rightarrow
1$, and
\item[(2)] root-$n$ consistency: $\hat{\mathbf {B}}_l = \mathbf{
B}_l^*\mathbf {D}_l^0 + \mathrm{O}_P(n^{-1/2})$ for $l = 1, \ldots, g$.
\end{longlist}
\end{theo}
%

Theorem~\ref{theo1ch5}, part (1), demonstrates that the shrinkage
group-wise MAVE method can efficiently remove unimportant predictors,
while part (2) implies that the estimator that corresponds to relevant
predictors is root-$n$ consistent. As we can see, the result is very
similar to that of adaptive LASSO for linear models (Zou
\cite{Zou06}). In
fact, we shall show in the proof that shrinkage group-wise MAVE is
closely related to an adaptive LASSO problem. A similar phenomena can
be found in Bondell and Li \cite{BonLi09} where they
studied the shrinkage
inverse regression estimation.
However, unlike linear
models, we need to study the interplay between the bandwidth $h$ and
the penalty parameter $\lambda$. This is explicitly shown in Theorem~\ref{theo1ch5}
 in which we require that $\lambda\rightarrow\infty$ and $\lambda
n^{-1/2}h^{-2} \rightarrow0$.
We also note that, although it is
possible to derive the asymptotic distribution, the form of the
asymptotic variance is rather complicated and thus is not pursued here.

As a direct application we consider the special case when $g = 1$, that
is, there is no group information available. It follows that the
shrinkage MAVE estimator possesses exactly the same properties.

%
%
\begin{coll}
Assume that $g = 1$, and that the regularity conditions \textup{(A1)--(A6)} given
in the \hyperref[app]{Appendix} hold. If $\lambda_n \rightarrow\infty$ and $\lambda
_n n^{-1/2} h^{-2} \rightarrow0$, then we have
\begin{longlist}[(2)]
\item[(1)] selection consistency: $P\{\mathcal{I}(\hat{\mathbf {B}})
= \mathcal{I}(\mathbf{ B}^*)\} \rightarrow1$, and
\item[(2)] root-$n$ consistency: $\hat{\mathbf {B}} = \mathbf{
B}^*\mathbf{ D}^0 + \mathrm{O}_P(n^{-1/2})$.
\end{longlist}
\end{coll}
%

The attractive properties of shrinkage group-wise MAVE depend
critically on an appropriate choice of the tuning parameter, for
which prediction based criteria such as generalized cross-validation
have been commonly used in practice. However, it is well known that
this practice tends to produce over-fitted models. For model
selection consistency, it has been verified that tuning parameter
selectors with the Bayesian information criterion are able to
identify the true model consistently; see for example Wang, Li and Tsai
 \cite{WanLiTsa07} and Wang and Leng
\cite{WanLen07}. In the following, we propose a
criterion which is similar in form to the classical Bayesian
information criterion.

Let $\hat{\alp}(\lambda) = (\hat{\alp}_1^{\top}, \ldots,
\hat{\alp}_g^{\top})^{\top}$. Write $\hat{\alp}(\lambda) =
(\hat{\alpha}_1, \ldots, \hat{\alpha}_p)^{\top}$ and $\bfitV= (V_1,
\ldots, V_p)^{\top}$. We use the notation $\mathcal{M} = \{r_{1},
r_{2}, \ldots, r_{p^*}\} \subseteq\{1, \ldots, p\}$ to denote an
arbitrary candidate model which includes predictors $\{V_{s}, s \in
\mathcal{M}\}$. Let $k_1 = 0$ and $k_l = p_1 + \cdots+ p_{l - 1}$
for $l = 2, \ldots, g$. Then, $\mathcal{M}_F = \{1, \ldots, p\}$ and
$\mathcal{M}_T = \bigcup_{l = 1}^g\{k_l + 1, \ldots, k_l + q_l\}$
represent the full model and the true model, respectively. Finally,
we use $|\mathcal{M}|$ to denote the size of the model
$\mathcal{M}$.

Let $\mathcal{M}_{\lambda} = \{s: \hat{\alpha}_s \neq0\}$ be the
model that is identified by $\hat{\alp}(\lambda)$ or $\hat{\B}$.
Define
\[
\operatorname{RSS}_{\lambda} = \sum_{i = 1}^{n}
\sum_{j = 1}^{n} \Biggl\{y^j - {
\tilde{a}^i} - \sum_{l = 1}^{g}{
\tilde{\bfitb}}_l^{i\top} \tilde{\tilde{\B}}_l{\vphantom{\tilde{B}}\!\!}^{\top}
\operatorname{diag}\bigl(\boldintav _l^j - \boldintav
_l^i\bigr)\hat{\alp}_l \Biggr
\}^2 \tilde{w}_{j}^i.
\]
We select the optimal $\lambda$ by minimizing
%
%
\begin{equation}
\label{BICch5} \operatorname{BIC}_{\lambda} = \log(\operatorname{RSS}_{\lambda})
+ \operatorname{df}_{\lambda} \frac{\log(n)}{n},
\end{equation}
where $\operatorname{df}_{\lambda}$ denotes the effective number of parameters
in the shrinkage group-wise MAVE estimator. The resulting optimal
regularization parameter is denoted by $\hat{\lambda}_{\mathrm{BIC}}$.
Following the discussion of Zou, Hastie and Tibshirani \cite{ZouHasTib07} about the
degrees of freedom of the LASSO estimator, we approximate $\operatorname{df}_{\lambda}$ by $d_1|\mathcal{M}_{\lambda}^1| + \cdots+
d_g|\mathcal{M}_{\lambda}^g|$, where $\mathcal{M}_{\lambda}^l$
represents the index set of identified predictors in group $l$.

We now establish the asymptotic property of the BIC-type tuning
parameter selector.

%
\begin{theo}\label{theo2ch5}
Suppose that the regularity conditions \textup{(A1)--(A6)} given
in the \hyperref[app]{Appendix}
hold. Then we have $P(\mathcal{M}_{\hat{\lambda}_{\mathrm{BIC}}} =
\mathcal{M}_T) \rightarrow1$.
\end{theo}

\begin{Remarks*}
\begin{enumerate}[3.]
\item Mixed-rates behavior naturally arises in the estimation of
semi-parametric models. As shown in the proof of Theorem~\ref
{theo1ch5}, the objective function (\ref{2.2.2ch5}) can be decomposed
into two components with different convergent rates. As a result, the
standard approach does not yield the complete limiting behavior of the
estimator. Fortunately, we are able to derive the asymptotic behavior
by directly applying results from mixed-rates asymptotics (Radchenko
\cite{Rad08}).

\item In practice, one may use a concave penalty other than the LASSO
penalty. We have tried using the smoothly clipped absolute deviation
penalty (Fan and Li \cite{FanLi01}) and the minimax
concave penalty (Zhang \cite{Zha10}),
and have found that the resulting estimators enjoy the same properties.
See Section~\ref{sec3ch5} for a numerical comparison of these methods.
Consider again the illustrative example in Section~\ref{sec1ch5}, the
proposed sparse group-wise MAVE method, when the smoothly clipped
absolute deviation penalty is used, yielded the direction estimates
\[
(0.702, -0.712, 0, \ldots, 0)^{\top} \quad\mbox{and}\quad (0.722, 0.692,
0, 0.027, 0, \ldots, 0)^{\top}.
\]
As we can see, all except one of the coordinates corresponding to
irrelevant predictors were correctly shrunk to zero.

\item The result here is applicable to a general class of
semi-parametric models. In particular, it provides an alternative
method for estimation and selection for partially linear single-index
models in which two groups exist naturally (Xia and H\"{a}rdle \cite{XiaHar06}).
Further, the new method can be adjusted to handle dimension reduction
and variable selection with censored data (Xia, Zhang and Xu
\cite{XiaZhaXu10}).
\item Although in this paper we focus on shrinkage estimation of the
group-wise central mean subspace, the same strategy can be used to
target the group-wise central subspace. To see this, we note that
Wang and Xia \cite{WanXia08} have modified MAVE to
estimate the central
subspace, and
so group-wise MAVE can be modified in a similar way to estimate the
group-wise central subspace; see Section~8 of Li, Li and Zhu
\cite{LiLiZhu10} for
more discussion. To conclude, we believe that these efforts would
enhance the usefulness of the shrinkage MAVE method in data analysis.
\end{enumerate}
\end{Remarks*}

\section{Numerical studies}\label{sec3ch5}

\subsection{Simulation studies}

In this subsection, we use simulations to evaluate the finite-sample
performance of the shrinkage group-wise MAVE method. For comparison we
consider the LASSO penalty, the smoothly clipped absolute deviation
(SCAD) penalty and the minimax concave penalty (MCP) in the simulation.
The resulting estimators, including group-wise MAVE, are denoted by
SgMAVE-LASSO, SgMAVE-SCAD, SgMAVE-MCP and gMAVE, respectively.
Throughout the following numerical studies we adopt the Gaussian kernel
and use the optimal bandwidth $h = \{4/(d+2)\}^{1/(d+4)}n^{-1/(d+4)}$.
The R code that we used for group-wise MAVE is available at
\url{http://www4.stat.ncsu.edu/\textasciitilde li/software.html}.
SgMAVE-LASSO is computed
using the least angle regression algorithm (Efron \textit{et al.} \cite{Efretal04}), while
SgMAVE-SCAD and SgMAVE-MCP are computed using the coordinate descent
algorithms described by Breheny and Huang \cite{BreHua11}.
The entire R code can
be requested from the authors.

To evaluate estimation accuracy, we compute the vector correlation
coefficient (VCC), which is defined as $(\prod_{t = 1}^{d_l} \phi
_t^2)^{1/ 2}$, and the trace correlation coefficient (TCC), which is
defined as $({d_l}^{-1}\sum_{t = 1}^{d_l} \phi_t^2)^{1 / 2}$, where
the $\phi_t^2$'s are the eigenvalues of the matrix $\hat{\B}_l^{\top
}\B_{l}^*\B_{l}^{*\top} \hat{\B}_l$. These two measures range
between 0 and 1, with larger values indicating a more
accurate estimator; see Ye and Weiss \cite{YeWei03} for
more information. We
also employ three summary statistics to assess how well the methods
select predictors: the average model size (MS), which is the average
number of nonzero rows of $\hat{\B}_l$; the true positive rate (TPR),
which is the average fraction of nonzero rows of $\hat{\B}_l$
associated with relevant predictors; and the false positive rate (FPR),
which is the average fraction of nonzero rows of $\hat{\B}_l$
associated with irrelevant predictors. Both TPR and FPR range between 0
and 1, and ideally, we wish to have TPR to be close to 1 and FPR to be
close to 0 at the same time. We report the results using the BIC-type
criterion (\ref{BICch5}) to select tuning parameters.

The predictor vector $\bfitV= (\bfitV_1^{\top}, \ldots, \bfitV
_g^{\top})^{\top}$ is generated from $N(\mathbf{0}_{p},
\boldsymbol{\Sigma})$ in each example. We examine two commonly-used
correlation structures among the predictors. The first is
autoregressive, $\Sigma_{st} = 0.5^{|s - t|}$ for all $s, t = 1,
\ldots, p$. Consequently, the predictors with large distances in order
are expected to be mutually independent approximately. The second is
compound symmetry, $\Sigma_{ss} = 1$ and $\Sigma_{st} = 0.5$ for any
$s \neq t$, so all the predictors are equally correlated with each other.

%
\begin{exmp}\label{ex1ch5} In this experiment, we set
\[
(g, p_1, p_2, d_1, d_2,
q_1, q_2) = (2, 20, 20, 1, 1, 3, 2).
\]
Thus, there are two groups, $\bfitV_1$ and $\bfitV_2$, and each group
consists of twenty predictors. Further, each predictor group is
connected with the response variable through a single linear
combination. Specifically, the response variable is generated from each
of the following three models:
%
%
\begin{eqnarray}
\label{mod11ch5}Y &=& \be_1^{\top}\bfitV_1\bigl(1
+ \be_2^{\top}\bfitV_2\bigr) + 0.5\epsilon,
\\
\label{mod12ch5}Y &=& \be_1^{\top}\bfitV_1/\bigl
\{0.5 + \bigl(1.5 + \be_2^{\top}\bfitV _2
\bigr)^2\bigr\} + 0.5\epsilon,
\\
\label{mod13ch5}Y &=& \exp\bigl(0.5 \be_1^{\top}
\bfitV_1\bigr) + \sin\bigl(0.2 \uppi\be_2^{\top
}
\bfitV_2\bigr) + 0.5\epsilon,
\end{eqnarray}
where $\be_1=(1, 1, 1, 0, \ldots, 0)^{\top}$, $\be_2=(1, 1, 0, \ldots
, 0)^{\top}$, $\epsilon\thicksim N(0, 1)$, and $\epsilon$ is
independent of all predictors. We let $n = 200$.
\end{exmp}

Table~\ref{tab1ch5} presents the simulation results based on 200 data
replications for these three models. As we can see, all methods
considered show very good performance, but the shrinkage ones often
achieve higher estimation accuracy than the one without shrinkage.
Further, although none of the three shrinkage methods can universally
dominate the other two competitors, SgMAVE-SCAD and SgMAVE-MCP tend to
produce sparser solutions than SgMAVE-LASSO. Finally, the performance
of the group-wise MAVE estimator and its shrinkage versions is only
slightly affected by the correlation structure among the predictors.

%
\begin{sidewaystable}
\tablewidth=\textwidth
\tabcolsep=0pt
\caption{Summary of Example \protect\ref{ex1ch5}. The
average vector correlation coefficient (VCC) with standard error in
parentheses, the average number of predictors selected (MS), true
positive rate (TPR) and false positive rate (FPR), based on 200 data
replications, are reported}\label{tab1ch5}
\begin{tabular*}{\textwidth}{@{\extracolsep{\fill}}lllllllll@{}}
\hline
&\multicolumn{4}{l}{$\be_{1}$}&\multicolumn{4}{l@{}}{$\be_{2}$}
\\[-5pt]
&\multicolumn{4}{l}{\hrulefill}&\multicolumn{4}{l@{}}{\hrulefill}
\\
& \multicolumn{1}{l}{VCC} & \multicolumn{1}{l}{MS}
& \multicolumn{1}{l}{TPR} & \multicolumn{1}{l}{FPR}
& \multicolumn{1}{l}{VCC} & \multicolumn{1}{l}{MS}
& \multicolumn{1}{l}{TPR} & \multicolumn{1}{l@{}}{FPR}
\\
\hline
& \multicolumn{8}{c@{}}{Model (\ref{mod11ch5}): autoregressive
correlation}
\\
gMAVE & 0.9963 (0.0239)& & & & 0.9905 (0.0780)& & &
\\
SgMAVE-LASSO& 0.9895 (0.0996)& 4.1550& 0.9900& 0.1692& 0.9896 (0.0997)&
3.1850& 0.9900& 0.1506
\\
SgMAVE-SCAD & 0.9976 (0.0207)& 3.0450& 1.0000& 0.0064& 0.9898 (0.0997)&
2.0150& 0.9900& 0.0043
\\
SgMAVE-MCP & 0.9977 (0.0193)& 3.1200& 1.0000& 0.0171& 0.9897 (0.0997)&
2.0550& 0.9900& 0.0093
\\[3pt]
& \multicolumn{8}{c@{}}{Model (\ref{mod11ch5}): compound symmetry}
\\
gMAVE & 0.9923 (0.0704)& & & & 0.9818 (0.1108)& & &
\\
SgMAVE-LASSO& 0.9933 (0.0709)& 4.8450& 0.9950& 0.2657& 0.9808 (0.1254)&
4.0300& 0.9900& 0.2562
\\
SgMAVE-SCAD & 0.9935 (0.0710)& 3.1400& 0.9950& 0.0221& 0.9811 (0.1252)&
2.1300& 0.9875& 0.0193
\\
SgMAVE-MCP & 0.9934 (0.0710)& 3.1950& 0.9950& 0.0300& 0.9805 (0.1294)&
2.1050& 0.9825& 0.0175
\\[5pt]
& \multicolumn{8}{c@{}}{Model (\ref{mod12ch5}): autoregressive
correlation}
\\
gMAVE & 0.9771 (0.0137)& & & & 0.9735 (0.0538)& & &
\\
SgMAVE-LASSO& 0.9885 (0.0104)& 5.5350& 1.0000& 0.3621& 0.9846 (0.0477)&
4.1400& 0.9975& 0.2681
\\
SgMAVE-SCAD & 0.9915 (0.0103)& 3.7300& 1.0000& 0.1042& 0.9849 (0.0557)&
2.5700& 0.9950& 0.0725
\\
SgMAVE-MCP & 0.9886 (0.0116)& 4.0100& 1.0000& 0.1442& 0.9837 (0.0550)&
2.5750& 0.9950& 0.0731
\\[3pt]
& \multicolumn{8}{c@{}}{Model (\ref{mod12ch5}): compound symmetry}
\\
gMAVE & 0.9739 (0.0177)& & & & 0.9432 (0.1535)& & &
\\
SgMAVE-LASSO& 0.9856 (0.0120)& 5.7650& 1.0000& 0.3950& 0.9450 (0.1940)&
3.7450& 0.9625& 0.2275
\\
SgMAVE-SCAD & 0.9896 (0.0130)& 3.8300& 1.0000& 0.1185& 0.9486 (0.1923)&
2.3500& 0.9600& 0.0537
\\
SgMAVE-MCP & 0.9858 (0.0132)& 4.0500& 1.0000& 0.1500& 0.9438 (0.2031)&
2.3000& 0.9550& 0.0487
\\
\hline
\end{tabular*}
\end{sidewaystable}

\setcounter{table}{0}
\begin{sidewaystable}
\tablewidth=\textwidth
\tabcolsep=0pt
\caption{(Continued)}
\begin{tabular*}{\textwidth}{@{\extracolsep{\fill}}lllllllll@{}}
\hline
&\multicolumn{4}{l}{$\be_{1}$}&\multicolumn{4}{l@{}}{$\be_{2}$}
\\[-5pt]
&\multicolumn{4}{l}{\hrulefill}&\multicolumn{4}{l@{}}{\hrulefill}
\\
& \multicolumn{1}{l}{VCC} & \multicolumn{1}{l}{MS}
& \multicolumn{1}{l}{TPR} & \multicolumn{1}{l}{FPR}
& \multicolumn{1}{l}{VCC} & \multicolumn{1}{l}{MS}
& \multicolumn{1}{l}{TPR} & \multicolumn{1}{l@{}}{FPR}
\\
\hline
& \multicolumn{8}{c@{}}{Model (\ref{mod13ch5}): autoregressive
correlation}
\\
gMAVE & 0.9955 (0.0026)& & & & 0.9648 (0.0214)& & &
\\
SgMAVE-LASSO& 0.9981 (0.0019)& 5.2600& 1.0000& 0.3228& 0.9879 (0.0158)&
3.6250& 1.0000& 0.2031
\\
SgMAVE-SCAD & 0.9984 (0.0019)& 3.8250& 1.0000& 0.1178& 0.9832 (0.0726)&
2.6850& 0.9950& 0.0868
\\
SgMAVE-MCP & 0.9981 (0.0020)& 3.6450& 1.0000& 0.0921& 0.9874 (0.0191)&
2.5650& 1.0000& 0.0706
\\[3pt]
& \multicolumn{8}{c@{}}{Model (\ref{mod13ch5}): compound symmetry}
\\
gMAVE & 0.9954 (0.0023)& & & & 0.9546 (0.0401)& & &
\\
SgMAVE-LASSO& 0.9974 (0.0019)& 5.8050& 1.0000& 0.4007& 0.9766 (0.0384)&
3.9250& 1.0000& 0.2406
\\
SgMAVE-SCAD & 0.9981 (0.0019)& 4.0600& 1.0000& 0.1514& 0.9792 (0.0517)&
2.8800& 0.9975& 0.1106
\\
SgMAVE-MCP & 0.9977 (0.0021)& 3.7850& 1.0000& 0.1121& 0.9786 (0.0499)&
2.5950& 0.9975& 0.0750
\\
\hline
\end{tabular*}
\end{sidewaystable}

%
\begin{exmp}\label{ex2ch5} In this experiment, we set
\[
(g, p_1, p_2, d_1, d_2,
q_1, q_2) = (2, 20, 20, 2, 1, q_1, 2).
\]
Thus, there are two groups, $\bfitV_1$ and $\bfitV_2$, and each group
consists of twenty predictors. Further, the first predictor group is
connected with the response variable through two linear combinations
and the second predictor group is connected with the response variable
through a single linear combination. The regression model is
%
%
\begin{equation}
\label{mod2ch5} Y = 2.5 \be_{11}^{\top}\bfitV_1/
\bigl\{0.5 + \bigl(1.5 + \be_{12}^{\top
}\bfitV_1
\bigr)^2\bigr\} + \be_2^{\top}\bfitV_2
+ 0.5\epsilon,
\end{equation}
where $\be_{11}=(1, 1, 0, \ldots, 0)^{\top}$, $\be_2=(1, 1, 0, \ldots
, 0)^{\top}$, $\epsilon\thicksim N(0, 1)$, and $\epsilon$ is
independent of all predictors. We consider two cases. In Case 1: we set
$q_1 = 2$ and $\be_{12}=(1, -1, 0, \ldots, 0)^{\top}$. In Case 2: we
set $q_1 = 4$ and $\be_{12}=(0, 0, 1, 1, 0, \ldots, 0)^{\top}$. We
let $n = 200$.
\end{exmp}

Table~\ref{tab2ch5} summarizes the simulation results out of 200 data
replications for Case 1 and Case 2. As in the previous example, we have
the same observations. Unreported results also show that the BIC-type
criterion has a pretty large rate of correctly identifying the true
model in this example.

%
\begin{sidewaystable}
\tablewidth=\textwidth
\tabcolsep=0pt
\caption{Summary of Example \protect\ref{ex2ch5}. The
average of the vector correlation coefficient (VCC) and the trace
correlation coefficient (TCC) with standard errors in parentheses, the
average number of predictors selected (MS), true positive rate (TPR)
and false positive rate (FPR), based on 200 data replications, are
reported}\label{tab2ch5}
\begin{tabular*}{\textwidth}{@{\extracolsep{\fill}}llllllllll@{}}
\hline
&\multicolumn{5}{l}{$\be_{1} = (\be_{11}, \be_{12})$}&\multicolumn
{4}{l@{}}{$\be_{2}$}
\\[-5pt]
&\multicolumn{5}{l}{\hrulefill}&\multicolumn{4}{l@{}}{\hrulefill}
\\
& \multicolumn{1}{l}{VCC} & \multicolumn{1}{l}{TCC} & \multicolumn{1}{l}{MS}
& \multicolumn{1}{l}{TPR} & \multicolumn{1}{l}{FPR}
& \multicolumn{1}{l}{VCC} & \multicolumn{1}{l}{MS}
& \multicolumn{1}{l}{TPR} & \multicolumn{1}{l@{}}{FPR}
\\
\hline
& \multicolumn{9}{c@{}}{Model (\ref{mod2ch5}): Case 1,
autoregressive correlation}
\\
gMAVE & 0.9506 (0.0972)& 0.9762 (0.0397)& & & & 0.9667 (0.0249)& & &
\\
SgMAVE-LASSO& 0.9936 (0.0707)& 0.9978 (0.0207)& 3.4150& 1.0000& 0.1271&
0.9867 (0.0997)& 2.2450& 0.9900& 0.0147
\\
SgMAVE-SCAD & 0.9916 (0.0772)& 0.9968 (0.0256)& 2.4350& 1.0000& 0.0755&
0.9917 (0.0734)& 2.0900& 0.9925& 0.0058
\\
SgMAVE-MCP & 0.9897 (0.0774)& 0.9973 (0.0158)& 2.6900& 1.0000& 0.0889&
0.9912 (0.0734)& 2.1700& 0.9925& 0.0102
\\[3pt]
& \multicolumn{9}{c@{}}{Model (\ref{mod2ch5}): Case 1, compound
symmetry}
\\
gMAVE & 0.9515 (0.0709)& 0.9760 (0.0312)& & & & 0.9616 (0.0191)& & &
\\
SgMAVE-LASSO& 0.9950 (0.0170)& 0.9975 (0.0085)& 4.4250& 1.0000& 0.1802&
0.9933 (0.0224)& 2.9200& 0.9975& 0.0513
\\
SgMAVE-SCAD & 0.9953 (0.0195)& 0.9976 (0.0097)& 2.9250& 1.0000& 0.1013&
0.9958 (0.0213)& 2.1700& 0.9975& 0.0097
\\
SgMAVE-MCP & 0.9879 (0.0302)& 0.9940 (0.0149)& 3.4700& 1.0000& 0.1300&
0.9933 (0.0136)& 2.5150& 1.0000& 0.0286
\\[5pt]
& \multicolumn{9}{c@{}}{Model (\ref{mod2ch5}): Case 2,
autoregressive correlation}
\\
gMAVE & 0.9548 (0.0764)& 0.9775 (0.0348)& & & & 0.9686 (0.0176)& & &
\\
SgMAVE-LASSO& 0.9498 (0.1843)& 0.9786 (0.0877)& 8.2150& 0.9787& 0.2687&
0.9887 (0.0739)& 4.6850& 0.9925& 0.1500
\\
SgMAVE-SCAD & 0.9723 (0.1289)& 0.9891 (0.0428)& 6.2700& 0.9925& 0.1437&
0.9952 (0.0183)& 2.7600& 1.0000& 0.0422
\\
SgMAVE-MCP & 0.9797 (0.0837)& 0.9907 (0.0319)& 5.6850& 0.9975& 0.1059&
0.9939 (0.0165)& 2.6150& 1.0000& 0.0341
\\[3pt]
& \multicolumn{9}{c@{}}{Model (\ref{mod2ch5}): Case 2, compound
symmetry}
\\
gMAVE & 0.9584 (0.0559)& 0.9795 (0.0219)& & & & 0.9648 (0.0167)& & &
\\
SgMAVE-LASSO& 0.9827 (0.0724)& 0.9923 (0.0228)& 8.8900& 0.9975& 0.3062&
0.9918 (0.0162)& 5.1550& 1.0000& 0.1752
\\
SgMAVE-SCAD & 0.9871 (0.0722)& 0.9945 (0.0224)& 5.3550& 0.9975& 0.0853&
0.9962 (0.0108)& 2.5300& 1.0000& 0.0294
\\
SgMAVE-MCP & 0.9831 (0.0727)& 0.9925 (0.0230)& 5.4800& 0.9975& 0.0931&
0.9934 (0.0146)& 2.5200& 1.0000& 0.0288
\\
\hline
\end{tabular*}
\end{sidewaystable}

%
\begin{exmp}\label{ex3ch5} In this experiment, we set
\[
(g, p_1, p_2, p_3, d_1,
d_2, d_3, q_1, q_2,
q_3) = (3, p_0, p_0, p_0, 1, 1,
1, 2, 2, 2).
\]
Thus, there are three groups, $\bfitV_1$, $\bfitV_2$ and $\bfitV_3$,
and each group consists of $p_0$ predictors. Further, each predictor
group is connected with the response variable through a single linear
combination. We consider the following two models:
%
%
\begin{eqnarray}
\label{mod31ch5}Y &=& \be_{1}^{\top}\bfitV_1 + 2
\be_{2}^{\top}\bfitV_2 / \bigl\{0.5 + \bigl(1.5 +
\be_{3}^{\top}\bfitV_3\bigr)^2\bigr\}
+ 0.5\epsilon,
\\
\label{mod32ch5}Y &=& \be_{1}^{\top}\bfitV_1 + 0.2
\bigl(2 + \be_{2}^{\top}\bfitV_2
\bigr)^2 + 2 \sin\bigl(0.2 \uppi\be_{3}^{\top}
\bfitV_3\bigr) + 0.5\epsilon,
\end{eqnarray}
where $\be_{1}=(1, -1, 0, \ldots, 0)^{\top}$, $\be_{2}=(1, 1, 0,
\ldots, 0)^{\top}$, $\be_3=(1, -1, 0, \ldots, 0)^{\top}$, $\epsilon
\thicksim N(0, 1)$, and $\epsilon$ is independent of all predictors.
We let $(n, p_0)$ be $(100, 10)$, $(200, 20)$ and $(200, 30)$.
\end{exmp}

The simulation results for models (\ref{mod31ch5}) and
(\ref{mod32ch5}), based on the 200 data replications, are shown in
Tables~\ref{tab3ch5} and~\ref{tab4ch5}, respectively. In general, the
results show that reasonably, increasing the sample size improves the
performance, while increasing the dimension of predictors makes the
performance worse. Moreover, the empirical performance of the shrinkage
estimators relies on the initial estimator as expected. Thus, the
development of a shrinkage estimation and variable selection method
that depends less on the initial estimator can be practically useful,
and we will work along this line in our future study.

\begin{sidewaystable}
\tablewidth=\textwidth
\tabcolsep=0pt
\caption{Summary of Example \protect\ref{ex3ch5}. The
average vector correlation coefficient (VCC) with standard error in
parentheses, the average number of predictors selected (MS), true
positive rate (TPR) and false positive rate (FPR), based on 200 data
replications, for model (\protect\ref{mod31ch5}), are reported}\label{tab3ch5}
\begin{tabular*}{\textwidth}{@{\extracolsep{\fill}}lllllllllllll@{}}
\hline
&\multicolumn{4}{l}{$\be_{1}$}&\multicolumn{4}{l}{$\be
_{2}$}&\multicolumn{4}{l@{}}{$\be_{3}$}
\\[-5pt]
&\multicolumn{4}{l}{\hrulefill}&\multicolumn{4}{l}{\hrulefill}
&\multicolumn{4}{l@{}}{\hrulefill}
\\
& \multicolumn{1}{l}{VCC} & \multicolumn{1}{l}{MS}
& \multicolumn{1}{l}{TPR} & \multicolumn{1}{l}{FPR}
& \multicolumn{1}{l}{VCC} & \multicolumn{1}{l}{MS}
& \multicolumn{1}{l}{TPR} & \multicolumn{1}{l}{FPR}
& \multicolumn{1}{l}{VCC} & \multicolumn{1}{l}{MS}
& \multicolumn{1}{l}{TPR} & \multicolumn{1}{l@{}}{FPR}
\\
\hline
& \multicolumn{12}{c@{}}{Model (\ref{mod31ch5}): $(n, p_0) =
(100, 10)$, autoregressive correlation}
\\
gMAVE & 0.966 (0.021)& & & & 0.965 (0.023)& & & & 0.956 (0.100)& & &
\\
SgMAVE-LASSO& 0.967 (0.106)& 5.850& 0.980& 0.486& 0.980 (0.025)& 7.060&
0.997& 0.633& 0.965 (0.127)& 6.075& 0.980& 0.514
\\
SgMAVE-SCAD & 0.952 (0.170)& 4.230& 0.967& 0.286& 0.980 (0.021)& 4.745&
1.000& 0.343& 0.960 (0.147)& 4.225& 0.977& 0.283
\\
SgMAVE-MCP & 0.972 (0.102)& 3.585& 0.987& 0.201& 0.979 (0.021)& 4.180&
1.000& 0.272& 0.960 (0.146)& 3.675& 0.980& 0.214
\\[3pt]
& \multicolumn{12}{c@{}}{Model (\ref{mod31ch5}): $(n, p_0) = (100, 10)$,
compound symmetry}
\\
gMAVE & 0.949 (0.046)& & & & 0.950 (0.049)& & & & 0.888 (0.225)& & &
\\
SgMAVE-LASSO& 0.935 (0.169)& 5.135& 0.950& 0.404& 0.953 (0.128)& 6.270&
0.975& 0.540& 0.881 (0.282)& 5.245& 0.915& 0.426
\\
SgMAVE-SCAD & 0.952 (0.130)& 3.950& 0.977& 0.249& 0.967 (0.052)& 4.560&
1.000& 0.320& 0.892 (0.258)& 3.875& 0.932& 0.251
\\
SgMAVE-MCP & 0.966 (0.052)& 3.460& 0.995& 0.183& 0.961 (0.081)& 4.100&
0.995& 0.263& 0.888 (0.264)& 3.420& 0.922& 0.196
\\[5pt]
 & \multicolumn{12}{c@{}}{Model (\ref{mod31ch5}): $(n, p_0) =
(200, 20)$, autoregressive correlation}
\\
gMAVE & 0.971 (0.018)& & & & 0.972 (0.014)& & & & 0.967 (0.087)& & &
\\
SgMAVE-LASSO& 0.950 (0.198)& 4.425& 0.952& 0.140& 0.994 (0.004)& 8.840&
1.000& 0.380& 0.960 (0.184)& 5.645& 0.962& 0.206
\\
SgMAVE-SCAD & 0.972 (0.142)& 4.545& 0.975& 0.144& 0.992 (0.009)& 6.060&
1.000& 0.225& 0.979 (0.121)& 4.900& 0.985& 0.162
\\
SgMAVE-MCP & 0.987 (0.071)& 3.525& 0.995& 0.085& 0.989 (0.010)& 5.125&
1.000& 0.173& 0.978 (0.121)& 3.905& 0.985& 0.107
\\[3pt]
& \multicolumn{12}{c@{}}{Model (\ref{mod31ch5}): $(n, p_0) = (200, 20)$,
compound symmetry}
\\
gMAVE & 0.961 (0.031)& & & & 0.962 (0.028)& & & & 0.926 (0.175)& & &
\\
SgMAVE-LASSO& 0.926 (0.235)& 3.505& 0.925& 0.091& 0.978 (0.104)& 7.505&
0.985& 0.307& 0.912 (0.269)& 4.550& 0.915& 0.151
\\
SgMAVE-SCAD & 0.942 (0.213)& 3.625& 0.950& 0.095& 0.982 (0.073)& 5.080&
0.995& 0.171& 0.925 (0.241)& 3.990& 0.935& 0.117
\\
SgMAVE-MCP & 0.983 (0.029)& 3.535& 1.000& 0.085& 0.980 (0.026)& 4.855&
1.000& 0.158& 0.926 (0.232)& 3.690& 0.935& 0.101
\\
\hline
\end{tabular*}
\end{sidewaystable}

\setcounter{table}{2}
\begin{sidewaystable}
\tablewidth=\textwidth
\tabcolsep=0pt
\caption{(Continued)}
\begin{tabular*}{\textwidth}{@{\extracolsep{\fill}}lllllllllllll@{}}
\hline
&\multicolumn{4}{l}{$\be_{1}$}&\multicolumn{4}{l}{$\be
_{2}$}&\multicolumn{4}{l@{}}{$\be_{3}$}
\\[-5pt]
&\multicolumn{4}{l}{\hrulefill}&\multicolumn{4}{l}{\hrulefill}
&\multicolumn{4}{l@{}}{\hrulefill}
\\
& \multicolumn{1}{l}{VCC} & \multicolumn{1}{l}{MS}
& \multicolumn{1}{l}{TPR} & \multicolumn{1}{l}{FPR}
& \multicolumn{1}{l}{VCC} & \multicolumn{1}{l}{MS}
& \multicolumn{1}{l}{TPR} & \multicolumn{1}{l}{FPR}
& \multicolumn{1}{l}{VCC} & \multicolumn{1}{l}{MS}
& \multicolumn{1}{l}{TPR} & \multicolumn{1}{l@{}}{FPR}
\\
\hline
& \multicolumn{12}{c@{}}{Model (\ref{mod31ch5}): $(n, p_0) =
(200, 30)$, autoregressive correlation}
\\
gMAVE & 0.943 (0.026)& & & & 0.938 (0.030)& & & & 0.943 (0.038)& & &
\\
SgMAVE-LASSO& 0.793 (0.375)& 2.895& 0.787& 0.047& 0.993 (0.005)& 8.500&
1.000& 0.232& 0.874 (0.308)& 4.200& 0.867& 0.088
\\
SgMAVE-SCAD & 0.913 (0.257)& 4.485& 0.912& 0.095& 0.985 (0.021)& 7.535&
1.000& 0.197& 0.977 (0.121)& 5.285& 0.985& 0.118
\\
SgMAVE-MCP & 0.976 (0.121)& 3.710& 0.985& 0.062& 0.979 (0.022)& 6.790&
1.000& 0.171& 0.990 (0.014)& 4.395& 1.000& 0.085
\\[3pt]
& \multicolumn{12}{c@{}}{Model (\ref{mod31ch5}): $(n, p_0) = (200, 30)$,
compound symmetry}
\\
gMAVE & 0.907 (0.058)& & & & 0.901 (0.063)& & & & 0.836 (0.223)& & &
\\
SgMAVE-LASSO& 0.828 (0.333)& 2.350& 0.817& 0.025& 0.952 (0.173)& 6.535&
0.960& 0.164& 0.827 (0.353)& 3.295& 0.822& 0.058
\\
SgMAVE-SCAD & 0.848 (0.328)& 2.930& 0.845& 0.044& 0.963 (0.125)& 5.250&
0.982& 0.117& 0.837 (0.347)& 3.750& 0.845& 0.073
\\
SgMAVE-MCP & 0.943 (0.172)& 3.460& 0.965& 0.054& 0.960 (0.039)& 6.240&
1.000& 0.151& 0.867 (0.298)& 3.920& 0.887& 0.076
\\
\hline
\end{tabular*}
\end{sidewaystable}

%

%
\begin{sidewaystable}
\tablewidth=\textwidth
\tabcolsep=0pt
\caption{Summary of Example \protect\ref{ex3ch5}. The
average vector correlation coefficient (VCC) with standard error in
parentheses, the average number of predictors selected (MS), true
positive rate (TPR) and false positive rate (FPR), based on 200 data
replications, for model (\protect\ref{mod32ch5}), are reported}\label{tab4ch5}
 \begin{tabular*}{\textwidth}{@{\extracolsep{\fill}}lllllllllllll@{}}
\hline
&\multicolumn{4}{l}{$\be_{1}$}&\multicolumn{4}{l}{$\be
_{2}$}&\multicolumn{4}{l@{}}{$\be_{3}$}
\\[-5pt]
&\multicolumn{4}{l}{\hrulefill}&\multicolumn{4}{l}{\hrulefill}
&\multicolumn{4}{l@{}}{\hrulefill}
\\
& \multicolumn{1}{l}{VCC} & \multicolumn{1}{l}{MS}
& \multicolumn{1}{l}{TPR} & \multicolumn{1}{l}{FPR}
& \multicolumn{1}{l}{VCC} & \multicolumn{1}{l}{MS}
& \multicolumn{1}{l}{TPR} & \multicolumn{1}{l}{FPR}
& \multicolumn{1}{l}{VCC} & \multicolumn{1}{l}{MS}
& \multicolumn{1}{l}{TPR} & \multicolumn{1}{l@{}}{FPR}
\\
\hline
& \multicolumn{12}{c@{}}{Model (\ref{mod32ch5}): $(n, p_0) =
(100, 10)$, autoregressive correlation}
\\
gMAVE & 0.968 (0.088)& & & & 0.965 (0.069)& & & & 0.979 (0.026)& & &
\\
SgMAVE-LASSO& 0.962 (0.154)& 6.970& 0.975& 0.627& 0.977 (0.071)& 7.410&
0.995& 0.677& 0.968 (0.140)& 7.105& 0.977& 0.643
\\
SgMAVE-SCAD & 0.973 (0.120)& 4.570& 0.987& 0.324& 0.981 (0.028)& 4.685&
0.997& 0.336& 0.970 (0.139)& 4.685& 0.980& 0.340
\\
SgMAVE-MCP & 0.973 (0.121)& 3.880& 0.985& 0.238& 0.978 (0.071)& 4.060&
0.995& 0.258& 0.970 (0.139)& 4.070& 0.980& 0.263
\\[3pt]
& \multicolumn{12}{c@{}}{Model (\ref{mod32ch5}): $(n, p_0) = (100, 10)$,
compound symmetry}
\\
gMAVE & 0.975 (0.014)& & & & 0.966 (0.052)& & & & 0.978 (0.010)& & &
\\
SgMAVE-LASSO& 0.985 (0.014)& 6.535& 1.000& 0.566& 0.978 (0.053)& 6.815&
0.997& 0.602& 0.987 (0.010)& 6.825& 1.000& 0.603
\\
SgMAVE-SCAD & 0.987 (0.014)& 4.030& 1.000& 0.253& 0.980 (0.071)& 3.945&
0.995& 0.244& 0.988 (0.011)& 4.320& 1.000& 0.290
\\
SgMAVE-MCP & 0.984 (0.013)& 3.845& 1.000& 0.230& 0.977 (0.071)& 3.835&
0.995& 0.230& 0.985 (0.010)& 4.150& 1.000& 0.268
\\[5pt]
 & \multicolumn{12}{c@{}}{Model (\ref{mod32ch5}): $(n, p_0) =
(200, 20)$, autoregressive correlation}
\\
gMAVE & 0.979 (0.007)& & & & 0.977 (0.009)& & & & 0.983 (0.006)& & &
\\
SgMAVE-LASSO& 0.997 (0.003)& 5.815& 1.000& 0.211& 0.995 (0.004)& 8.215&
1.000& 0.345& 0.997 (0.003)& 5.935& 1.000& 0.218
\\
SgMAVE-SCAD & 0.996 (0.006)& 4.505& 1.000& 0.139& 0.995 (0.006)& 4.705&
1.000& 0.150& 0.997 (0.005)& 4.570& 1.000& 0.142
\\
SgMAVE-MCP & 0.994 (0.006)& 4.240& 1.000& 0.124& 0.993 (0.006)& 4.230&
1.000& 0.123& 0.995 (0.005)& 4.445& 1.000& 0.135
\\[3pt]
& \multicolumn{12}{c@{}}{Model (\ref{mod32ch5}): $(n, p_0) = (200, 20)$,
compound symmetry}
\\
gMAVE & 0.976 (0.008)& & & & 0.974 (0.011)& & & & 0.980 (0.007)& & &
\\
SgMAVE-LASSO& 0.997 (0.002)& 3.845& 1.000& 0.102& 0.995 (0.004)& 6.340&
1.000& 0.241& 0.997 (0.003)& 4.315& 1.000& 0.128
\\
SgMAVE-SCAD & 0.996 (0.006)& 3.170& 1.000& 0.065& 0.996 (0.005)& 3.370&
1.000& 0.076& 0.997 (0.005)& 3.330& 1.000& 0.073
\\
SgMAVE-MCP & 0.991 (0.007)& 4.095& 1.000& 0.116& 0.989 (0.009)& 4.210&
1.000& 0.122& 0.992 (0.006)& 4.310& 1.000& 0.128
\\
\hline
\end{tabular*}
\end{sidewaystable}

\setcounter{table}{3}
\begin{sidewaystable}
\tablewidth=\textwidth
\tabcolsep=0pt
\caption{(Continued)}
\begin{tabular*}{\textwidth}{@{\extracolsep{\fill}}lllllllllllll@{}}
\hline
&\multicolumn{4}{l}{$\be_{1}$}&\multicolumn{4}{l}{$\be
_{2}$}&\multicolumn{4}{l@{}}{$\be_{3}$}
\\[-5pt]
&\multicolumn{4}{l}{\hrulefill}&\multicolumn{4}{l}{\hrulefill}
&\multicolumn{4}{l@{}}{\hrulefill}
\\
& \multicolumn{1}{l}{VCC} & \multicolumn{1}{l}{MS}
& \multicolumn{1}{l}{TPR} & \multicolumn{1}{l}{FPR}
& \multicolumn{1}{l}{VCC} & \multicolumn{1}{l}{MS}
& \multicolumn{1}{l}{TPR} & \multicolumn{1}{l}{FPR}
& \multicolumn{1}{l}{VCC} & \multicolumn{1}{l}{MS}
& \multicolumn{1}{l}{TPR} & \multicolumn{1}{l@{}}{FPR}
\\
\hline
 & \multicolumn{12}{c@{}}{Model (\ref{mod32ch5}): $(n, p_0) =
 (200, 30)$, autoregressive correlation}
\\
gMAVE & 0.960 (0.013)& & & & 0.951 (0.017)& & & & 0.967 (0.012)& & &
\\
SgMAVE-LASSO& 0.993 (0.070)& 3.680& 0.995& 0.060& 0.996 (0.003)& 7.470&
1.000& 0.195& 0.993 (0.070)& 3.840& 0.995& 0.066
\\
SgMAVE-SCAD & 0.997 (0.006)& 3.710& 1.000& 0.061& 0.996 (0.008)& 4.335&
1.000& 0.083& 0.998 (0.004)& 4.005& 1.000& 0.071
\\
SgMAVE-MCP & 0.995 (0.007)& 3.810& 1.000& 0.064& 0.992 (0.009)& 4.255&
1.000& 0.080& 0.995 (0.006)& 4.030& 1.000& 0.072
\\[3pt]
& \multicolumn{12}{c@{}}{Model (\ref{mod32ch5}): $(n, p_0) = (200, 30)$,
compound symmetry}
\\
gMAVE & 0.951 (0.014)& & & & 0.943 (0.021)& & & & 0.959 (0.015)& & &
\\
SgMAVE-LASSO& 0.997 (0.003)& 2.895& 1.000& 0.031& 0.995 (0.003)& 6.340&
1.000& 0.155& 0.998 (0.003)& 3.010& 1.000& 0.036
\\
SgMAVE-SCAD & 0.998 (0.003)& 2.345& 1.000& 0.012& 0.998 (0.003)& 2.635&
1.000& 0.022& 0.998 (0.002)& 2.375& 1.000& 0.013
\\
SgMAVE-MCP & 0.989 (0.009)& 4.190& 1.000& 0.078& 0.986 (0.012)& 4.820&
1.000& 0.100& 0.991 (0.007)& 4.310& 1.000& 0.082
\\
\hline
\end{tabular*}
\end{sidewaystable}

As we mentioned before, the MAVE procedure is a
combination of nonparametric function estimation and direction
estimation; it is an iterative procedure with each cycle consisting of
two least squares problems. As we know, inverse regression based
methods, which are largely ``parametric'' in nature, are simple and easy
to use. Thus, the proposed approach is computationally more demanding
than inverse regression based methods, especially when the sample size
and the predictor dimension are very high. Table~\ref{tab6} shows the
average CPU times, based on 200 data replications, for the shrinkage
group-wise MAVE method (along with penalty parameter selection) for
model (\ref{mod32ch5}) in Example~\ref{ex3ch5}. All algorithms are implemented as R language
functions, and all timings were carried out on a Dell Poweredge R410
dual processors server equipped with Six Core Xeon X5670 2.93~GHz CPU,
64~GB RAM running CentOS 5 Linux. We see that the times depend on both
$n$ and $p$. We also find similar results (unreported) for the other
models considered in the simulation studies. Nevertheless, we emphasize
that, as opposed to inverse regression based methods which require
strong conditions on the distribution of predictors, direct regression
based methods such as MAVE need relatively weak conditions such as the
smoothness of the link function, and they often have much better
performance for finite samples.

\begin{table}
\tablewidth=\textwidth
\tabcolsep=0pt
\caption{Run times (CPU seconds) for shrinkage group-wise MAVE of
various sizes $n$, $p$ and different correlation structures among the
predictors for model (\ref{mod32ch5}) in Example~\protect\ref{ex3ch5}}\label{tab6}
\begin{tabular*}{\textwidth}{@{\extracolsep{\fill}}lll@{}}
\hline
&Autoregressive correlation & Compound symmetry
\\
\hline
$(n, p) = (100, 30)$ & \phantom{0}14 & \phantom{0}13
\\
$(n, p) = (200, 60)$ & \phantom{0}87 & \phantom{0}99
\\
$(n, p) = (200, 90)$ & 186 & 200
\\
\hline
\end{tabular*}
\end{table}


\subsection{Pyrimidine data}

A common step in drug design is the formation of a quantitative
structure-activity relationship (QSAR; So \cite{So00}). The QSAR analysis is
to relate a numerical description of molecular structure to known
biological activity. The pyrimidine data set, which is available in the
UCI machine-learning repository at
\url{http://archive.ics.uci.edu/ml/machine-learning-databases/qsar/}, was
studied by Hirst, King and Sternberg \cite{HirKinSte94} to model
the QSAR of the
inhibition of dihydrofolate reductase (DHFR) by pyrimidines. It
contains a structural information on 74 2,4-diamino-5-(substituted
benzyl)pyrimidines used as inhibitors of DHFR in Escherichia coli. Each
pyrimidine compound has 3 positions of substitution where chemical
activity occurs, and at each position the substituent is assigned nine\vadjust{\goodbreak}
physicochemical attributes: polarity (PL), size (SZ), flexibility (FL),
number of hydrogen-bond donors (HD), number of hydrogen-bond acceptors
(HA), strength and presence of $\pi$-donors ($\pi$D), strength and
presence of $\pi$-acceptors ($\pi$A), polarisability of the molecular
orbitals (PO) and $\sigma$-effect ($\sigma$E). The response variable
is the experimentally assayed activity of the inhibitors.

\begin{figure}

\includegraphics{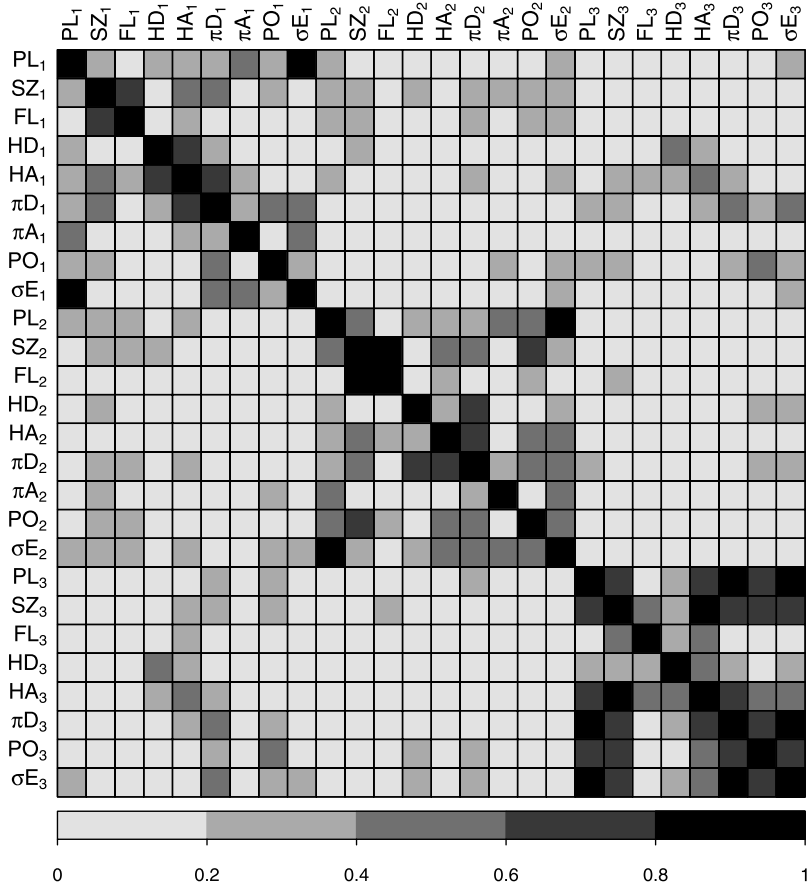}

\caption{Graphical representation of the absolute correlation matrix
of the $26$ predictors for the pyrimidine data. The magnitude of each
pairwise correlation is represented by a block in the grayscale
image.}\label{fig1}
\end{figure}

The attributes in this data set naturally fall into 3 groups
corresponding to the substitution positions 1, 2 and 3. This is further
confirmed by the graphical representation of the correlation matrix of
the attributes in Figure~\ref{fig1}; for example, the attributes
belonging to the third substitution position have a moderately strong
correlation, but are weakly associated with most of the other
attributes. We write $\bfitV_l = (\mathrm{PL}_l, \mathrm{SZ}_l,
\mathrm{FL}_l, \mathrm{HD}_l, \mathrm{HA}_l, \pi\mathrm{D}_l, \pi
\mathrm{A}_l, \mathrm{PO}_l, \sigma\mathrm{E}_l)^{\top}$ for $l =
1, 2$ and 3, where the attributes are represented by their two-letter
abbreviations with the subscripts denoting the position of
substitution. All predictors are standardized to have mean zero and
unit length (the predictor $\pi\mathrm{A}_3$ has no variability and
is then removed). Thus, in this data set, the sample size is $n = 74$,
the predictor dimension is $p = 26$, and the group information is $(g,
p_1, p_2, p_3) = (3, 9, 9, 8)$. We regard this ``prior'' information on
predictor group structure as a given fact.

Ordinary least squares (OLS), LASSO, SCAD, MCP, group-wise MAVE (gMAVE)
and shrinkage group-wise MAVE (SgMAVE-LASSO, SgMAVE-SCAD and
SgMAVE-MCP) are applied to this data set. Before applying the
group-wise MAVE procedure, we need to determine $(d_1, d_2, d_3)$. The
BIC-type criterion of Li, Li and Zhu \cite{LiLiZhu10}
that is a modification of
Wang and Yin \cite{WanYin08} yields $(d_1, d_2, d_3) =
(1, 1, 1)$, indicating
that each predictor group is connected with the response variable
through a single linear combination. The same criterion, when the group
information is ignored ($g = 1$), also shows a three-dimensional
structure in regression.

The corresponding coefficient estimates are shown in the second through
ninth columns of Table~\ref{tab5}. Using the attribute representation
(that is, representing molecules by a set of physicochemical
attributes), these methods generate a variety of possible influences of
structure on activity. As we can see, the substituent at position 1
should not be a hydrogen-bond acceptor ($\mathrm{HA}_1$). We can also
see that both the size and the flexibility of the substituent at
positions 1 and 3, say $\mathrm{SZ}_1, \mathrm{FL}_1, \mathrm{SZ}_3$
and $ \mathrm{FL}_3$, are informative to the activity of the
pyrimidines. Previous analysis of the crystal structure of the complex
formed between trimethoprim and DHFR shows that the substituents at
positions 1 and 3, are buried in a hydrophobic environment, and
restrictions on size and flexibility are consistent with this
(Hirst, King and Sternberg \cite{HirKinSte94}). The shrinkage
group-wise MAVE methods also
identify the $\sigma$-effect of the substituent at position 1 ($\sigma
\mathrm{E}_1$) as an influencing factor, which is in accordance with
previous studies using machine learning techniques. However, the
ordinary variable selection methods fail to detect it.

%


Let $\hat{\B}_1 \in\mathbb{R}^9$, $\hat{\B}_2 \in\mathbb{R}^9$
and $\hat{\B}_3 \in\mathbb{R}^8$ denote direction estimates for the
three predictor groups, respectively. We next consider the group-wise
additive index model
\[
Y = G_1(Z_1) + G_2(Z_2) +
G_3(Z_3) + \varepsilon,
\]
%
%
\begin{sidewaystable}
\tablewidth=\textwidth
\tabcolsep=0pt
\caption{Pyrimidine data. Estimated coefficients and adjusted
$R$-squared values ($\bar{R}^2$) from various methods}\label{tab5}
\begin{tabular*}{\textwidth}{@{\extracolsep{\fill}}ld{2.4}d{2.4}d{2.4}d{2.4}d{2.4}d{2.4}d{2.4}d{2.4}@{}}
\hline
&\multicolumn{1}{l}{OLS}&\multicolumn{1}{l}{LASSO}&\multicolumn
{1}{l}{SCAD}&\multicolumn{1}{l}{MCP}&\multicolumn
{1}{l}{gMAVE}&\multicolumn{1}{l}{SgMAVE-LASSO}
& \multicolumn{1}{l}{SgMAVE-SCAD} & \multicolumn{1}{l@{}}{SgMAVE-MCP}
\\
\hline
& \multicolumn{8}{c@{}}{$\bfitV_1$: attributes of a substituent at
position 1}
\\
$\mathrm{PL}_1$ & -0.0437& 0 & 0 & 0 & 0.6240& 0.5620& 0.6043& 0.5896
\\
$\mathrm{SZ}_1$ & 0.0435& 0.0358& 0 & 0.0586& -0.3154& -0.3772&
-0.3336& -0.3470
\\
$\mathrm{FL}_1$ & -0.0444& -0.0323& -0.0216& -0.0537& 0.2890& 0.3444&
0.3035& 0.3198
\\
$\mathrm{HD}_1$ & -0.0249& -0.0197& -0.0480& -0.0273& 0.1079& 0.1589&
0.1288& 0.1407
\\
$\mathrm{HA}_1$ & 0.0139& 0 & 0.0358& 0 & -0.0227& 0 & 0 & 0
\\
$\pi\mathrm{D}_1$ & 0.0083& 0.0025& 0.0185& 0.0207& -0.0521& -0.0928&
-0.0898& -0.0916
\\
$\pi\mathrm{A}_1$ & 0.0122& 0 & 0.0160& 0 & -0.1901& -0.1818&
-0.1985& -0.1949
\\
$\mathrm{PO}_1$ & 0.0288& 0.0224& 0.0291& 0 & -0.3498& -0.3604&
-0.3701& -0.3711
\\
$\sigma\mathrm{E}_1$ & 0.0386& 0.0080& 0 & 0 & -0.5040& -0.4756&
-0.4798& -0.4749
\\[5pt]
& \multicolumn{8}{c@{}}{$\bfitV_2$: attributes of a substituent at
position 2}
\\
$\mathrm{PL}_2$ & -0.0281& 0 & -0.0025& 0 & 0.6924& 0.6763& 0.6894&
0.6812
\\
$\mathrm{SZ}_2$ & 0.0396& 0.0152& 0.0324& 0 & -0.1058& -0.1023&
-0.0991& -0.0849
\\
$\mathrm{FL}_2$ & -0.0430& -0.0240& -0.0407& -0.0179& -0.1802&
-0.2211& -0.1904& -0.2232
\\
$\mathrm{HD}_2$ & 0.0076& 0 & 0 & 0 & -0.0982& -0.0556& -0.0841&
-0.0574
\\
$\mathrm{HA}_2$ & 0.0068& 0 & 0 & 0 & 0.0402& 0.0440& 0.0328& 0.0494
\\
$\pi\mathrm{D}_2$ & -0.0132& 0.0100& 0 & 0 & 0.0663& 0 & 0.0534& 0
\\
$\pi\mathrm{A}_2$ & 0.0030& 0 & 0 & 0 & -0.1654& -0.1702& -0.1672&
-0.1770
\\
$\mathrm{PO}_2$ & 0.0218& 0.0145& 0.0148& 0.0257& -0.3264& -0.3575&
-0.3348& -0.3462
\\
$\sigma\mathrm{E}_2$ & 0.0176& 0 & 0 & 0 & -0.5720& -0.5667& -0.5722&
-0.5673
\\[5pt]
& \multicolumn{8}{c@{}}{$\bfitV_3$: attributes of a substituent at
position 3}
\\
$\mathrm{PL}_3$ & 0.0103& 0 & 0 & 0 & -0.0925& -0.1599& -0.1506&
-0.1501
\\
$\mathrm{SZ}_3$ & 0.0458& 0.0579& 0.0955& 0.0740& 0.7768& 0.8369&
0.9125& 0.9126
\\
$\mathrm{FL}_3$ & -0.0311& -0.0079& -0.0309& -0.0208& -0.1757&
-0.1957& -0.2048& -0.2070
\\
$\mathrm{HD}_3$ & -0.0540& -0.0424& -0.0294& -0.0340& 0.0586& 0 & 0 &
0
\\
$\mathrm{HA}_3$ & 0.0595& 0 & 0 & 0 & -0.2928& -0.2298& -0.2689&
-0.2667
\\
$\pi\mathrm{D}_3$ & -0.1958& -0.0696& -0.1620& -0.1459& -0.3503&
-0.3490& -0.1738& -0.1744
\\
$\mathrm{PO}_3$ & 0.0696& 0.0365& 0.0470& 0.0623& -0.0992& 0.0080& 0 &
0
\\
$\sigma\mathrm{E}_3$ & 0.0987& 0.0363& 0.0864& 0.0806& 0.3677&
0.2468& 0 & 0
\\[3pt]
$\bar{R}^2$ & 0.8206& 0.7934& 0.8311& 0.8226& 0.9150& 0.9241& 0.9170&
0.9210
\\
\hline
\end{tabular*}%
\end{sidewaystable}
where $Z_l = \hat{\B}_l^{\top}\bfitV_l$, $l = 1, 2$ and $3$, are
the extracted linear predictors, and $G_l(\cdot)$'s are unknown
univariate functions. We fit this model by applying the \textbf{gam}
function in the publicly available R package \textbf{mgcv}. The
adjusted percentages of total deviance explained, namely the adjusted
\mbox{R-squared} values, for various methods are summarized in the last row of
Table~\ref{tab5}. Unreported results show that the nonparametric
smoothing of all the three predictors yields better performance than
the additive model using smoothing of every single predictor, but the
improvement is not statistically significant. As we can see, the
proposed semi-parametric methods outperform the classical parametric
ones as they can provide a mechanism for exploring nonlinear
relationships between molecular structure and biological activity.

\begin{figure}\vspace*{5pt}

\includegraphics{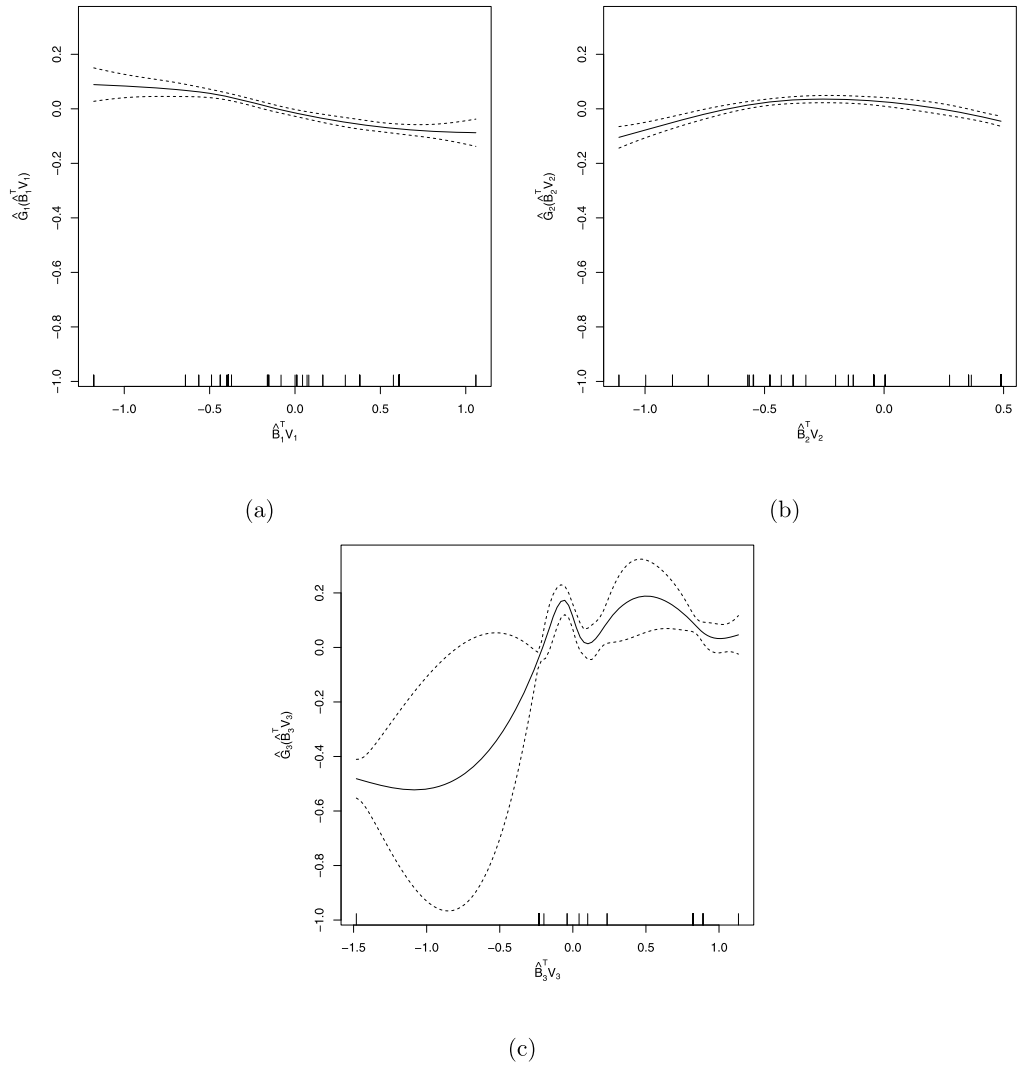}

\caption{The panels show the estimates of the terms in the group-wise
additive index model for the pyrimidine data using for illustration
SgMAVE-MCP. The upper left panel, the upper right panel and the lower
panel are the smooth functions of the extracted linear predictor in
predictor group one, two and three, respectively. The rug plots, along
the bottom of each plot, show the values of the predictors of each
smooth. Thin plate regression splines were used with smoothing
parameters being selected by GCV.}\label{fig2}
\end{figure}

Figure~\ref{fig2} provides the plots of estimated index functions,
using for illustration the shrinkage group-wise MAVE method with the
minimax concave penalty (SgMAVE-MCP). From Figures~\ref{fig2}(a) and~(b), it can be seen that $G_1(\cdot)$ has a linear trend, while
$G_2(\cdot)$ is clearly curved, indicating a nonlinear parabolic
dependence of activity on the extracted linear combination of
attributes at the second position of substitution. It can also be seen
from Figure~\ref{fig2}(c) that $G_3(\cdot)$ is very complicated, and
nonparametric smoothing performs poorly in areas where observations are sparse.

%

\section{Discussion}

In this paper, we only provide the convergence rate of
the shrinkage group-wise MAVE estimator. It is possible to derive the
limiting distribution. However, the limiting distribution is too
complicated to be applied for inference. Thus, for the time being, we
are frustrated by the lack of a good approximation to the limiting
distribution that can be used to set standard errors or to carry out
tests on the parameter vector.

As remarked by Knight and Fu \cite{KniFu00}, attaching
standard errors to
LASSO-type estimators is nontrivial. They then considered using the
residual-based bootstrap method to estimate the sampling distribution
of the LASSO estimator in a multiple linear regression. However,
Chatterjee and Lahiri \cite{ChaLah10} showed that the
conditional residual
bootstrap distribution given the data converges to a random measure;
that is, the residual bootstrap estimate of the LASSO distribution is
inconsistent. In a subsequent paper, Chatterjee and Lahiri \cite{ChaLah11}
proposed a modified bootstrap method, and showed that it provides a
valid approximation to the distribution of the LASSO estimator.

But it is unclear yet whether or not the modified bootstrap method of
Chatterjee and Lahiri \cite{ChaLah11} can be applied to our
setting. The
situation is complicated by the fact that in semi-parametric
multiple-index models we need to take into account the interaction
between nonparametric function estimation and shrinkage direction
estimation. Work along this line is in progress.

\begin{appendix}\label{app}
\section*{Appendix}

We need the following regularity conditions:
\begin{enumerate}[(A3)]
\item[(A1)] $E|Y|^k < \infty$ and $E\|\bfitX\|_2^k < \infty$ for
some large $k > 0$, where $\|\cdot\|_2$ denotes the $\ell_2$ norm.
\item[(A2)] The density function of $\bfitX$ has a bounded second
derivative; $E(\bfitX|\B^{\top}\bfitX= \bfitw)$ and $E(\bfitX
\bfitX^{\top}|\B^{\top}\bfitX= \bfitw)$ have bounded derivatives
with respect to $\bfitw$ and $\B$ for $\B$ in a small neighborhood
of $\B^*$, that is, $\|\P_{\B} - \P_{\B^*}\|_2
\leq\zeta$ for some small $\zeta> 0$.
\item[(A3)] The function $E(Y|\B^{\top}\bfitX= \bfitw)$ has a
bounded and continuous fourth derivative with respect to $\bfitw$ and
$\B$ for $\B$ in a small neighborhood of $\B^*$.
\item[(A4)] The kernel $K(\cdot)$ is a Gaussian probability density function.
\item[(A5)] $d \leq3$ and $h \propto n^{-1/(d + 4)}$.
\item[(A6)] $\tilde{\B}_l = \B_l^*\D_l^0 + \mathrm{O}_P(n^{-1/2})$ and
$\tilde{\tilde{\B}}_l = \B_l^*\D_l^0 + \mathrm{O}_P(n^{-1/2})$ for some
$d_l \times d_l$ nonsingular matrix $\D_l^0$ for $l = 1, \ldots, g$.
\end{enumerate}
Note that conditions (A1)--(A4) are standard in the literature, see for
instance Wang and Xia \cite{WanXia08}, Xia \cite{Xia08} and Li, Li and Zhu \cite{LiLiZhu10}. As
shown in Xia \cite{Xia08}, the ordinary MAVE
estimator is root-$n$ consistent
under conditions similar to (A1)--(A5). Consequently, condition (A6) is
very reasonable because we can view $\tilde{\B} = \bigoplus_{l = 1}^g
\tilde{\B}_l$ as a special case of a general $\tilde{\B}$ in Xia's
proof. If higher order local polynomial smoothing is used, the root-$n$
consistency can also be achieved for $d > 3$; see Remark~5.3 in
Xia \cite{Xia08}. Nevertheless, in practice
models with $d > 3$ are not
attractive due to the ``curse of dimensionality''.

Before we begin the proof, we need to introduce some additional
notation. For a positive integer $m$, $\mathbf{0}_m$
stands for an $m$-dimensional vector of zeros. For an $m_1 \times m_2$
matrix $\A$, $\operatorname{vec}(\A)$ stands for the $m_1m_2$-dimensional
vector obtained by stacking the columns of $\A$. For a diagonal
matrix, we get the (generalized) inverse by taking the reciprocal of
each nonzero element on the diagonal, leaving the zeros in place, and
transposing the resulting matrix.

Let $\hat{\B}_l = \operatorname{diag}(\hat{\alp}_l)\tilde{\tilde{\B
}}_l$. Then $\hat{\B} = \bigoplus_{l = 1}^g \hat{\B}_l$. Let $\tilde
{\tilde{B}}_{lst}$ denote the $(s, t)$th element of $\tilde{\tilde
{\B}}_l$. Without loss of generality, we assume that $\D_l^0 = \I
_{d_l}$, and the first $q_l$ components of $\B_{l1}^*$ are nonzero.
For each $l = 1, \ldots, g$, we define
\[
\H_{l}^* = \pmatrix{ \operatorname{diag}\bigl(\B^*_{l1}
\bigr)
\cr
\vdots
\cr
\operatorname{diag}\bigl(\B^*_{ld_l}\bigr)} \bigl\{
\operatorname{diag}\bigl(\B^*_{l1}\bigr)\bigr\}^{-1}\quad \mbox{\textrm{and}}\quad
\tilde{\tilde{\H}}_{l} = \pmatrix{\operatorname{diag}(\tilde{\tilde{
\B}}_{l1})
\cr
\vdots
\cr
\operatorname{diag}(\tilde{\tilde{
\B}}_{ld_l})} \bigl\{\operatorname{diag}(\tilde{\tilde{
\B}}_{l1})\bigr\}^{-1},
\]
where $\B^*_{lt}$ denotes the $t$th column of $\B^*_l$ and $\tilde
{\tilde{\B}}_{lt}$ denotes the $t$th column of $\tilde{\tilde{\B
}}_l$, $t = 1, \ldots, d_l$. Let $\H^* = \bigoplus_{l = 1}^g \H_l^*$
and $\tilde{\tilde{\H}} = \bigoplus_{l = 1}^g \tilde{\tilde{\H
}}_l$. By condition (A6), $\|\H^* - \tilde{\tilde{\H}}\|_2 = \mathrm{O}_P(n^{-1/2})$.

\begin{pf*}{Proof of Theorem~\ref{theo1ch5}}
We shall
concentrate on the optimization problem (\ref{2.2.2ch5}). The proof
follows Theorem~1 of Bondell and Li \cite{BonLi09}
closely. First, we formulate
an equivalent optimization problem that is easier to analyze
theoretically. To see this, we note that
\begin{eqnarray*}
&& \sum_{i = 1}^{n}\sum
_{j = 1}^{n} \Biggl\{y^j - {
\tilde{a}^i} - \sum_{l = 1}^{g}{
\tilde{\bfitb}}_l^{i\top} \tilde{\tilde{\B
}}_l{\vphantom{\tilde{B}}\!\!}^{\top} \operatorname{diag}\bigl(\boldintav
_l^j - \boldintav _l^i\bigr)
\alp_l \Biggr\}^2 \tilde{w}_{j}^i
\\
&&\quad = \sum_{i = 1}^{n}\sum
_{j = 1}^{n} \Biggl[y^j - {
\tilde{a}^i} - \sum_{l = 1}^{g}
\bigl\{{\tilde{\bfitb}}_l^{i\top} \otimes\bigl(\boldintav
_l^j - \boldintav _l^i
\bigr)^{\top
}\bigr\} \pmatrix{\operatorname{diag}(\tilde{\tilde{
\B}}_{l1})
\cr
\vdots
\cr
\operatorname{diag}(\tilde{\tilde{
\B}}_{ld_l})} \alp_l \Biggr]^2
\tilde{w}_{j}^i
\\
&&\quad = \sum_{i = 1}^{n}\sum
_{j = 1}^{n} \Biggl[y^j - {
\tilde{a}^i} - \sum_{l = 1}^{g}
\bigl\{{\tilde{\bfitb}}_l^{i\top} \otimes\bigl(\boldintav
_l^j - \boldintav _l^i
\bigr)^{\top
}\bigr\} \tilde{\tilde{\H}}_l\operatorname{diag}(
\tilde{\tilde{\B}}_{l1})\alp _l \Biggr]^2
\tilde{w}_{j}^i.
\end{eqnarray*}
Suppose that $\{\check{\B}_{l1} = (\check{B}_{l11}, \ldots, \check
{B}_{lp_l 1})^{\top} \in\mathbb{R}^{p_l}, l = 1, \ldots, g\}$ is
the minimizer of
%
%
\begin{equation}
\label{app1ch5} \sum_{i = 1}^{n}\sum
_{j = 1}^{n} \Biggl[y^j - {
\tilde{a}^i} - \sum_{l = 1}^{g}
\bigl\{{\tilde{\bfitb}}_l^{i\top} \otimes\bigl(\boldintav
_l^j - \boldintav _l^i
\bigr)^{\top}\bigr\} \tilde{\tilde{\H}}_l \B_{l1}
\Biggr]^2 \tilde{w}_{j}^i + \lambda
_n \sum_{l = 1}^g\sum
_{s = 1}^{p_l}\frac{|B_{ls1}|}{|\tilde{\tilde
{B}}_{ls1}|}
\end{equation}
with respect to $\{\B_{l1} = (B_{l11}, \ldots, B_{lp_l 1})^{\top}
\in\mathbb{R}^{p_l}, l = 1, \ldots, g\}$. From Definition~\ref
{defi1ch5}, it is easy to see that $\hat{\alpha}_{ls} = (\tilde
{\tilde{B}}_{ls1})^{-1}\check{B}_{ls1}$ for all $l = 1, \ldots, g$
and $s = 1, \ldots, p_l$. Further, $\operatorname{vec}({\hat{\B}}_l) =
\tilde{\tilde{\H}}_l \check{\B}_{l1}$.

Below we shall describe the details of the proof by breaking it up into
two steps. Step\vspace*{2pt} I establishes the convergence rate of $\hat{\B}$.
Step II shows that $\hat{\B}$ attains sparsity.

\emph{Step} I. Let $\bfitu= (\bfitu_1^{\top}, \ldots, \bfitu
_g^{\top})^{\top} \in\mathbb{R}^p$, where $\bfitu_l = (u_{l1},
\ldots, u_{lp_l})^{\top} \in\mathbb{R}^{p_l}$ for $l = 1, \ldots,
g$. Define
\[
J_n(\bfitu) = \lambda_n \sum
_{l = 1}^g\sum_{s = 1}^{p_l}|
\tilde {\tilde{B}}_{ls1}|^{-1}\biggl\llvert B_{ls1}^*
+ \frac{u_{ls}}{\sqrt
{n}}\biggr\rrvert .
\]
Then, we have
\[
J_n(\bfitu) - J_n(\mathbf{0}_p) =
\frac{\lambda
_n}{\sqrt{n}} \sum_{l = 1}^g\sum
_{s = 1}^{p_l}|\tilde{\tilde {B}}_{ls1}|^{-1}
\sqrt{n} \biggl(\biggl\llvert B_{ls1}^* + \frac{u_{ls}}{\sqrt
{n}}\biggr\rrvert
- \bigl |B_{ls1}^*\bigr | \biggr).
\]
If $u_{ls} = 0$, then
\[
\frac{\lambda_n}{\sqrt{n}}|\tilde{\tilde{B}}_{ls1}|^{-1}\sqrt {n}
\biggl(\biggl\llvert B_{ls1}^* + \frac{u_{ls}}{\sqrt{n}}\biggr\rrvert -
\bigl |B_{ls1}^*\bigr | \biggr) = 0.
\]
If $u_{ls} \neq0$ and $B_{ls1}^* \neq0$, then $|\tilde{\tilde
{B}}_{ls1}|^{-1} \rightarrow_P |B_{ls1}^*|^{-1}$ and
\[
\sqrt{n} \biggl(\biggl\llvert B_{ls1}^* + \frac{u_{ls}}{\sqrt{n}}\biggr
\rrvert - \bigl |B_{ls1}^*\bigr | \biggr) \rightarrow u_{ls} \times
\operatorname{sgn}\bigl(B_{ls1}^*\bigr),
\]
where $\operatorname{sgn}(\cdot)$ is the sign function. By Slutsky's theorem,
\[
\frac{\lambda_n}{\sqrt{n}} |\tilde{\tilde{B}}_{ls1}|^{-1}\sqrt {n}
\biggl(\biggl\llvert B_{ls1}^* + \frac{u_{ls}}{\sqrt{n}}\biggr\rrvert -
\bigl |B_{ls1}^*\bigr | \biggr) = \mathrm{o}_P\bigl(h^2
\bigr).
\]
If $u_{ls} \neq0$ and $B_{ls1}^* = 0$, then
\[
\frac{\lambda_n}{\sqrt{n}}|\tilde{\tilde{B}}_{ls1}|^{-1}\sqrt {n}
\biggl(\biggl\llvert B_{ls1}^* + \frac{u_{ls}}{\sqrt{n}}\biggr\rrvert -
\bigl |B_{ls1}^*\bigr | \biggr) = \frac{\lambda_n}{\sqrt{n}|\tilde{\tilde
{B}}_{ls1}|} |u_{ls}|
\rightarrow_P \infty.
\]

Define
\[
\Psi_n(\bfitu) = \sum_{i = 1}^{n}
\sum_{j = 1}^{n} \Biggl[y^j - {
\tilde{a}^i} - \sum_{l = 1}^{g}
\bigl\{{\tilde{\bfitb}}_l^{i\top} \otimes\bigl(\boldintav
_l^j - \boldintav _l^i
\bigr)^{\top}\bigr\} \tilde{\tilde{\H}}_l \biggl(
\B_{l1}^* + \frac
{\bfitu_l}{\sqrt{n}} \biggr) \Biggr]^2
\tilde{w}_{j}^i.
\]
After some algebra one gets
\begin{eqnarray*}
&&\Psi_n(\bfitu) - \Psi_n(\mathbf{0}_{p})
\\[1pt]
&&\quad= \sum_{i =
1}^{n}\sum
_{j = 1}^{n} \Biggl[\sum_{l = 1}^{g}
\bigl\{{\tilde{\bfitb }}_l^{i\top} \otimes\bigl(\boldintav
_l^j - \boldintav _l^i
\bigr)^{\top}\bigr\} \tilde{\tilde{\H}}_l\frac
{\bfitu_l}{\sqrt{n}}
\Biggr]^2 \tilde{w}_{j}^i
\\[1pt]
&&\qquad{}+ 2\sum_{i = 1}^{n}\sum
_{j = 1}^{n} \Biggl[\sum_{l = 1}^{g}
\bigl\{{\tilde {\bfitb}}_l^{i\top} \otimes\bigl(\boldintav
_l^j - \boldintav _l^i
\bigr)^{\top}\bigr\} \tilde{\tilde{\H }}_l\frac{\bfitu_l}{\sqrt{n}}
\Biggr]
\\[1pt]
&&\phantom{\qquad{}+} {}\times \Biggl[\sum_{l =
1}^{g}
\bigl\{{\tilde{\bfitb}}_l^{i\top} \otimes\bigl(\boldintav
_l^j - \boldintav _l^i
\bigr)^{\top}\bigr\} \tilde{\tilde {\H}}_l\bigl(
\B_{l1}^* - \tilde{\tilde{\B}}_{l1}\bigr) \Biggr]
\tilde{w}_{j}^i.
\end{eqnarray*}
Let $\boldintav ^* = (\B_{11}^{*\top}, \ldots, \B
_{g1}^{*\top})^{\top}$ and $\tilde{\tilde{\boldintav }} = (\tilde{\tilde{\B}}_{11}{\vphantom{\tilde{B}}\!\!}^{\top}, \ldots, \tilde{\tilde
{\B}}_{g1}{\vphantom{\tilde{B}}\!\!}^{\top})^{\top}$. Then, we have
\begin{eqnarray*}
\Psi_n(\bfitu) - \Psi_n(\mathbf{0}_{p}) &=&
\bfitu ^{\top} \tilde{\tilde{\H}}{\vphantom{\tilde{H}}\!\!}^{\top} \Biggl(\frac{1}{n}
\sum_{i =
1}^{n}\sum
_{j = 1}^{n} \tilde{\boldintav }_{ij}\tilde{
\boldintav }_{ij}^{\top} \tilde {w}_{j}^i
\Biggr) \tilde{\tilde{\H}} \bfitu
\\[1pt]
&&{} + 2 \bfitu^{\top} \tilde{\tilde{\H}}{\vphantom{\tilde{H}}\!\!}^{\top} \Biggl(
\frac{1}{n}\sum_{i = 1}^{n}\sum
_{j = 1}^{n} \tilde{\boldintav
}_{ij}\tilde{\boldintav }_{ij}^{\top}
\tilde{w}_{j}^i \Biggr) \tilde{\tilde{\H}} \bigl\{\sqrt{n}
\bigl(\boldintav ^* - \tilde{\tilde{\boldintav }}\bigr)\bigr\}
\\[1pt]
&\equiv& T_1 + T_2,
\end{eqnarray*}
where $\tilde{\boldintav }_{ij} = (\tilde{\boldintav }_{ij1}^{\top}, \ldots, \tilde{\boldintav }_{ijg}^{\top})^{\top} \in\mathbb{R}^{p_1d_1 + \cdots
+ p_g d_g}$, $\tilde{\boldintav }_{ijl} = {\tilde
{\bfitb}_l^i} \otimes(\boldintav _l^j - \boldintav _l^i) \in\mathbb{R}^{p_ld_l}$, $l = 1, \ldots, g$.

First, we consider $T_1$. Note that
\begin{eqnarray*}
T_1 &=& \bfitu^{\top}\H^{*\top} \Biggl(
\frac{1}{n}\sum_{i =
1}^{n}\sum
_{j = 1}^{n} \tilde{\boldintav
}_{ij}\tilde{\boldintav }_{ij}^{\top} \tilde
{w}_{j}^i \Biggr)\H^*\bfitu
\\
& &{} + \bfitu^{\top}\H^{*\top} \Biggl(\frac{1}{n}\sum
_{i =
1}^{n}\sum_{j = 1}^{n}
\tilde{\boldintav }_{ij}\tilde{\boldintav }_{ij}^{\top}
\tilde {w}_{j}^i \Biggr) \bigl(\tilde{\tilde{\H}} - \H^*
\bigr)\bfitu
\\
& &{} + \bfitu^{\top}\bigl(\tilde{\tilde{\H}} - \H^*
\bigr)^{\top
} \Biggl(\frac{1}{n}\sum_{i = 1}^{n}
\sum_{j = 1}^{n} \tilde{\boldintav
}_{ij}\tilde{\boldintav }_{ij}^{\top}
\tilde{w}_{j}^i \Biggr)\H^*\bfitu
\\
& &{} + \bfitu^{\top}\bigl(\tilde{\tilde{\H}} - \H^*
\bigr)^{\top
} \Biggl(\frac{1}{n}\sum_{i = 1}^{n}
\sum_{j = 1}^{n} \tilde{\boldintav
}_{ij}\tilde{\boldintav }_{ij}^{\top}
\tilde{w}_{j}^i \Biggr) \bigl(\tilde{\tilde{\H}} - \H ^*
\bigr)\bfitu
\\
&\equiv& T_{11} + T_{12} + T_{13} +
T_{14}.
\end{eqnarray*}
By Lemma~4 in Wang and Xia \cite{WanXia08},
\[
\frac{1}{n}\sum_{i = 1}^{n}\sum
_{j = 1}^{n} \tilde{\boldintav
}_{ij}\tilde{\boldintav }_{ij}^{\top}
\tilde{w}_{j}^i = \boldsymbol{\Delta} +
\mathrm{o}_P(1),
\]
where $\boldsymbol{\Delta}$ is a nonnegative definite matrix. Hence, we
obtain
\[
T_{11} = \mathrm{O}_P(1),\qquad T_{12} =
\mathrm{O}_P\bigl(n^{-1/2}\bigr),\qquad T_{13} =
\mathrm{O}_P\bigl(n^{-1/2}\bigr) \quad\mbox{and}\quad
T_{14} = \mathrm{O}_P\bigl(n^{-1}\bigr).
\]
Next, we consider $T_2$. Note that $\tilde{\tilde{\boldintav }} = \boldintav ^* + \mathrm{O}_P(n^{-1/2})$ and
\begin{eqnarray*}
T_2 &=& 2 \bfitu^{\top} \H^{*\top} \Biggl(
\frac{1}{n}\sum_{i =
1}^{n}\sum
_{j = 1}^{n} \tilde{\boldintav
}_{ij}\tilde{\boldintav }_{ij}^{\top} \tilde
{w}_{j}^i \Biggr) \H^* \bigl\{\sqrt{n}\bigl(\boldintav ^*
- \tilde{\tilde{\boldintav }}\bigr)\bigr\}
\\
&&{} + 2 \bfitu^{\top} \H^* \Biggl(\frac{1}{n}\sum
_{i =
1}^{n}\sum_{j = 1}^{n}
\tilde{\boldintav }_{ij}\tilde{\boldintav }_{ij}^{\top}
\tilde {w}_{j}^i \Biggr) \bigl(\tilde{\tilde{\H}} - \H^*
\bigr) \bigl\{\sqrt{n}\bigl(\boldintav ^* - \tilde{\tilde{\boldintav }}\bigr)\bigr
\}
\\
&&{} + 2 \bfitu^{\top} \bigl(\tilde{\tilde{\H}} - \H^*
\bigr)^{\top} \Biggl(\frac{1}{n}\sum_{i = 1}^{n}
\sum_{j = 1}^{n} \tilde{\boldintav
}_{ij}\tilde{\boldintav }_{ij}^{\top}
\tilde{w}_{j}^i \Biggr) \H^* \bigl\{\sqrt{n}\bigl(
\boldintav ^* - \tilde{\tilde{\boldintav }}\bigr)\bigr\}
\\
&&{} + 2 \bfitu^{\top} \bigl(\tilde{\tilde{\H}} - \H^*
\bigr)^{\top} \Biggl(\frac{1}{n}\sum_{i = 1}^{n}
\sum_{j = 1}^{n} \tilde{\boldintav
}_{ij}\tilde{\boldintav }_{ij}^{\top}
\tilde{w}_{j}^i \Biggr) \bigl(\tilde{\tilde{\H}} - \H^*
\bigr) \bigl\{\sqrt{n}\bigl(\boldintav ^* - \tilde{\tilde {\boldintav }}\bigr)
\bigr\}
\\
&\equiv& T_{21} + T_{22} + T_{23} +
T_{24}.
\end{eqnarray*}
Thus, we arrive at
\[
T_{21} = \mathrm{O}_P(1),\qquad T_{22} =
\mathrm{O}_P\bigl(n^{-1/2}\bigr),\qquad T_{23} =
\mathrm{O}_P\bigl(n^{-1/2}\bigr) \quad\mbox{and}\quad
T_{24} = \mathrm{O}_P\bigl(n^{-1}\bigr).
\]

Let $L_n(\bfitu) = \Psi_n(\bfitu) + J_n(\bfitu)$. If $u_{ls} \neq
0$ for some $l \in\{1, \ldots, g\}$ and $s \in\{q_l + 1, \ldots,
p_l\}$, then $L_n(\bfitu) - L_n(\mathbf{0}_p) \rightarrow
_P \infty> 0$. So we assume in the sequel that $\bfitu\in\mathcal
{U}$, where
\[
\mathcal{U} = \bigl\{\bfitu\in\mathbb{R}^p: u_{ls} = 0
\mbox{ for all } l = 1, \ldots, g\mbox{ and }s = q_l + 1, \ldots,
p_l\bigr\}.
\]
It follows that $J_n(\bfitu) - J_n(\mathbf{0}_p) =
\mathrm{o}_P(h^2)$ for any $\bfitu\in\mathcal{U}$. Let $\bfitu\in\mathcal{U}$.

We consider the problem of minimizing $L_n(\bfitu)$ over $\mathcal
{U}$. Because $nh^4 \rightarrow\infty$, we obtain
\[
L_n(\bfitu) - L_n(\mathbf{0}_p) =
T_{11} + T_{21} + \mathrm{o}_P
\bigl(h^2\bigr).
\]
Let $(\B_{l}^*, \A_{l}^*)$ be an orthogonal matrix. Let $\C= \bigoplus_{l = 1}^g \{\I_{d_l} \otimes(\B_l^*, \A_{l}^*)\}$. Then, according
to Lemma~4 of Wang and Xia \cite{WanXia08}, the long
version, there exists a
$(\sum_{l = 1}^g d_l p_l) \times(\sum_{l = 1}^g d_l p_l)$
permutation matrix $\bfPi= \bigoplus_{l = 1}^g \bfPi_l$ such that
\begin{eqnarray*}
&&\frac{1}{n}\sum_{i = 1}^{n}\sum
_{j = 1}^{n} \tilde{\boldintav
}_{ij}\tilde{\boldintav }_{ij}^{\top}
\tilde{w}_{j}^i= \C\bfPi %
\pmatrix{ \boldsymbol{
\Delta}_{11n} & \boldsymbol{\Delta}_{12n}
\cr
\boldsymbol{
\Delta}_{21n} & \boldsymbol{\Delta}_{22n} } %
(\C
\bfPi)^{\top} + \mathrm{o}_P\bigl(h^2\bigr),
\end{eqnarray*}
where $h^{-2}\boldsymbol{\Delta}_{11n} \rightarrow_P \boldsymbol{\Delta}_{11}$, $h^{-2}\boldsymbol{\Delta}_{12n} \rightarrow_P \boldsymbol{\Delta}_{12}$,
$h^{-2}\boldsymbol{\Delta}_{21n} \rightarrow_P \boldsymbol{\Delta}_{21}$ and $\boldsymbol{\Delta}_{22n} \rightarrow_P \boldsymbol{\Delta}_{22}$.
Moreover, both $\boldsymbol{\Delta}_{11}$ and $\boldsymbol{\Delta}_{22}$ are positive definite.

Write $\C\bfPi= (\D_1, \D_2)$ with $\D_1$ being of order $(\sum_{l = 1}^g d_l p_l) \times(\sum_{l = 1}^g d_l^2)$. Let $\bfitz_1 =
\D_1^{\top}\H^*\bfitu$ and $\bfitz_2 = \D_2^{\top}\H^*\bfitu$.
Write $\bfitz= (\bfitz_1^{\top}, \bfitz_2^{\top})^{\top}$. Now
consider the function
\begin{eqnarray*}
G_n(\bfitz) &=& \bfitz^{\top} %
\pmatrix{
\boldsymbol{\Delta}_{11n} & \boldsymbol{\Delta}_{12n}
\cr
\boldsymbol{\Delta}_{21n} & \boldsymbol{\Delta}_{22n} }
\bfitz
\\
&&{}+ 2\bfitz^{\top} %
\pmatrix{ \boldsymbol{
\Delta}_{11n} & \boldsymbol{\Delta}_{12n}
\cr
\boldsymbol{
\Delta}_{21n} & \boldsymbol{\Delta}_{22n} } %
(\C
\bfPi)^{\top} \H^* \bigl\{\sqrt{n}\bigl(\boldintav ^* - \tilde{\tilde{
\boldintav }}\bigr)\bigr\} + \mathrm{o}_P\bigl(h^2\bigr).
\end{eqnarray*}
Denote $\check{\bfitz} = (\check{\bfitz}_1^{\top}, \check{\bfitz
}_2^{\top})^{\top}$ the minimizer of $G_n(\bfitz)$. It turns out
that the conditions of Theorem~1 of Radchenko \cite{Rad08} are satisfied and,
consequently, we have $\check{\bfitz}_1 = \mathrm{O}_P(1)$ and $\check{\bfitz
}_2 = \mathrm{O}_P(1)$. Over $\mathcal{U}$, because $\bfitu= (\H^{*\top}\H
^*)^{-1}\H^{*\top}(\D_1, \D_2)\bfitz$, we have $\check{\bfitu} =
\mathrm{O}_P(1)$. We thus conclude that there exists a minimizer $\{\check{\B
}_{l1}, l = 1, \ldots, g\}$ of (\ref{app1ch5}) such that $\|\check
{\B}_{l1} - \B_{l1}^*\|_2 = \mathrm{O}_P(n^{-1/2})$ for all $l = 1, \ldots, g$.

Since $\operatorname{vec}({{\B}}_l^*) = {\H}_l^* {\B}_{l1}^*$ and
$\operatorname{vec}({\hat{\B}}_l) = \tilde{\tilde{\H}}_l \check{\B}_{l1}$, by
triangular inequality we have
\[
\bigl \|\operatorname{vec}(\hat{\B}_l) - \operatorname{vec}\bigl(
\B_l^*\bigr)\bigr \|_2 \leq\bigl \|\operatorname{vec}\bigl({
\H}_l^* \check{\B}_{l1}\bigr) - \operatorname{vec}\bigl(
\H_l^* \B_{l1}^*\bigr)\bigr \| _2 + \bigl \|
\operatorname{vec}(\tilde{\tilde{\H}}_l \check{\B}_{l1}) -
\operatorname{vec}\bigl({\H}_l^* \check{\B}_{l1}\bigr)
\bigr \|_2.
\]
Therefore, $\|\hat{\B}_{l} - \B_{l}^*\|_2 = \mathrm{O}_P(n^{-1/2})$ for all
$l = 1, \ldots, g$.

\emph{Step} II. We show the variable selection consistency. Write
$\mathcal{A}_l = \mathcal{I}(\B^*_l)$ and $\mathcal{A}_{nl} =
\mathcal{I}(\hat{\B}_{l})$. For any $s \in\bigcup_{l = 1}^{g}\mathcal
{A}_l$, that is, $s \in\mathcal{A}_l$ for some $l$, the estimation
consistency result indicates that $\hat{\alpha}_{ls} \rightarrow_P
1$. Thus, $P(s \in\bigcup_{l = 1}^{g} \mathcal{A}_{nl}) \rightarrow1$.
It then suffices to show that for any $s' \notin\bigcup_{l =
1}^{g}\mathcal{A}_l$, $P(s' \in\bigcup_{l = 1}^{g}\mathcal{A}_{nl})
\rightarrow0$. Consider the event $\{s' \in\mathcal{A}_{nl}\}$. By
standard Karush--Kuhn--Tucker conditions for optimality, we know that\vspace*{-1pt}
\begin{eqnarray*}
&&\frac{2}{\sqrt{n}}\sum_{i = 1}^{n}\sum
_{j = 1}^{n} \Biggl[y^j - {
\tilde{a}^i} - \sum_{l = 1}^{g}
\bigl\{{\tilde{\bfitb}}_l^{i\top} \otimes\bigl(\boldintav
_l^j - \boldintav _l^i
\bigr)^{\top}\bigr\} \tilde{\tilde{\H}}_l\check{
\B}_{l1} \Biggr]
\\
&&\quad{} \times\bigl[\bigl\{{\tilde{\bfitb}}_l^{i\top}
\otimes\bigl(\boldintav _l^j - \boldintav
_l^i\bigr)^{\top}\bigr\} \tilde{\tilde{
\H}}_l \bfite_{ls'}\bigr] \times\tilde{w}_{j}^i
= |\tilde{\tilde{B}}_{ls'1}|^{-1}\operatorname{sgn}(
\check{B}_{ls'1})\frac
{\lambda_n}{\sqrt{n}},
\end{eqnarray*}
where $\bfite_{ls'} \in\mathbb{R}^{p_l}$ is the vector containing a
1 in the $s'$th position and zeros elsewhere. Note that\vspace*{-1pt}
\[
|\tilde{\tilde{B}}_{ls'1}|^{-1}\frac{\lambda_n}{\sqrt{n}} =
\frac
{\lambda_n}{\sqrt{n}|\tilde{\tilde{B}}_{ls'1}|} \rightarrow_P \infty
\]
and\vspace*{-1pt}
\begin{eqnarray*}
&&\frac{2}{\sqrt{n}}\sum_{i = 1}^{n}\sum
_{j = 1}^{n} \Biggl[y^j - {
\tilde{a}^i} - \sum_{l = 1}^{g}
\bigl\{{\tilde{\bfitb}}_l^{i\top} \otimes\bigl(\boldintav
_l^j - \boldintav _l^i
\bigr)^{\top}\bigr\} \tilde{\tilde{\H}}_l\check{
\B}_{l1} \Biggr]
\\
&&\quad{}\times\bigl[\bigl\{{\tilde{\bfitb}}_l^{i\top}
\otimes\bigl(\boldintav _l^j - \boldintav
_l^i\bigr)^{\top}\bigr\} \tilde{\tilde{
\H}}_l \bfite_{ls'}\bigr] \times\tilde{w}_{j}^i
= \mathrm{O}_P(1).
\end{eqnarray*}
%
Thus, we obtain\vspace*{-1pt}
\[
P\Biggl(s' \in\bigcup_{l = 1}^{g}
\mathcal{A}_{nl}\Biggr) \leq\sum_{l = 1}^{g}
P\bigl(s' \in\mathcal{A}_{nl}\bigr) \rightarrow0.
\]
The proof is complete.\vspace*{-1pt}
\end{pf*}

\begin{pf*}{Proof of Theorem~\ref{theo2ch5}}
According to
whether the fitted model $\mathcal{M}_{\lambda}$ is under-fitted,
correctly fitted or over-fitted, we can divide $\mathbb{R}^+ = [0,
\infty)$ into three disjoint parts:\vspace*{-1pt}
\[
\Omega_- = \{\lambda: \mathcal{M}_{\lambda} \nsupseteq
\mathcal{M}_T\},\qquad \Omega_0 = \{\lambda:
\mathcal{M}_{\lambda} = \mathcal{M}_T\}
\]
and\vspace*{-1pt}
\[
\Omega_+ = \{\lambda: \mathcal{M}_{\lambda} \supseteq\mathcal{M}_T,
\mathcal{M}_{\lambda} \neq\mathcal{M}_T\}.
\]
Further, we assume a reference sequence of tuning parameters, $\{
\lambda_n\}_{n = 1}^{\infty}$, which satisfies the conditions in
Theorem~\ref{theo1ch5}. Clearly, $\check{\B}_{l1} = \B_{l1}^* +
\mathrm{O}_P(n^{-1/2})$ and $P(\mathcal{M}_{\lambda_n} = \mathcal{M}_{T})
\rightarrow1$.

We write $\alp= (\alp_1^{\top}, \ldots, \alp_g^{\top})^{\top}
\equiv(\alpha_1, \ldots, \alpha_p)^{\top} \in\mathbb{R}^p$ and define\vspace*{-1pt}
\begin{eqnarray*}
\operatorname{RSS}_{\mathcal{M}} &=& \min_{\alp \in\mathcal
{S}_{\mathcal{M}}}\sum
_{i = 1}^{n}\sum_{j = 1}^{n}
\Biggl\{y^j - {\tilde{a}^i} - \sum
_{l = 1}^{g}{\tilde{\bfitb}}_l^{i\top}
\tilde {\tilde{\B}}_l{\vphantom{\tilde{B}}\!\!}^{\top} \operatorname{diag}\bigl(
\boldintav _l^j - \boldintav _l^i
\bigr)\alp_l \Biggr\}^2 \tilde {w}_{j}^i
\\
&\equiv& \min_{\alp \in\mathcal{S}_{\mathcal{M}}} \sum_{i = 1}^{n}
\sum_{j = 1}^{n} \Biggl[y^j - {
\tilde{a}^i} - \sum_{l =
1}^{g}
\bigl\{{\tilde{\bfitb}}_l^{i\top} \otimes\bigl(\boldintav
_l^j - \boldintav _l^i
\bigr)^{\top}\bigr\} \tilde{\tilde {\H}}_l
\operatorname{diag}(\tilde{\tilde{\B}}_{l1})\alp_l
\Biggr]^2 \tilde{w}_{j}^i,
\end{eqnarray*}
where $\mathcal{S}_{\mathcal{M}} = \{\bfitw= (w_1, \ldots,
w_p)^{\top} \in\mathbb{R}^p: w_s = 0, s \notin\mathcal{M}\}$.\vadjust{\goodbreak}

For a generic model $\mathcal{M}$, let $\{\breve{\B}_{11}(\mathcal
{M}), \ldots, \breve{\B}_{g1}(\mathcal{M})\}$ be the minimizer of
\[
\sum_{i = 1}^{n}\sum
_{j = 1}^{n} \Biggl[y^j - {
\tilde{a}^i} - \sum_{l = 1}^{g}
\bigl\{{\tilde{\bfitb}}_l^{i\top} \otimes\bigl(\boldintav
_l^j - \boldintav _l^i
\bigr)^{\top}\bigr\} \tilde{\tilde{\H}}_l \B_{l1}
\Biggr]^2 \tilde{w}_{j}^i
\]
with respect to $(\B_{11}^{\top}, \ldots, \B_{g1}^{\top})^{\top}
\in\mathcal{S}_{\mathcal{M}} $. Then, we have
\[
\operatorname{RSS}_{\mathcal{M}} = \sum_{i = 1}^{n}
\sum_{j = 1}^{n} \Biggl[y^j - {
\tilde{a}^i} - \sum_{l = 1}^{g}
\bigl\{{\tilde{\bfitb}}_l^{i\top
} \otimes\bigl(\boldintav
_l^j - \boldintav _l^i
\bigr)^{\top}\bigr\} \tilde{\tilde{\H}}_l \breve{
\B}_{l1} \Biggr]^2 \tilde{w}_{j}^i.
\]
Further, for the full model $\mathcal{M}_F$, $\breve{\B
}_{l1}(\mathcal{M}_F) = \tilde{\tilde{\B}}_{l1}$ for all $l = 1,
\ldots, g$.

We first consider under-fitted models, that is, $\mathcal{M}_{\lambda
} \nsupseteq\mathcal{M}_T$. Note that
\[
\inf_{\lambda\in\Omega_-} \operatorname{RSS}_{\lambda} -
\operatorname{RSS}_{\lambda_n} \geq\inf_{\lambda\in\Omega_-}
\operatorname{RSS}_{\mathcal{M}_{\lambda}} - \operatorname{RSS}_{\lambda_n} \geq\min
_{\mathcal{M} \nsupseteq\mathcal{M}_T} \operatorname{RSS}_{\mathcal
{M}} - \operatorname{RSS}_{\lambda_n}.
\]
By definition, we know that
\begin{eqnarray*}
\operatorname{RSS}_{\mathcal{M}} - \operatorname{RSS}_{\mathcal{M}_F} & = & \sum
_{i = 1}^{n}\sum_{j = 1}^{n}
\Biggl[y^j - {\tilde{a}^i} - \sum
_{l = 1}^{g}\bigl\{{\tilde{\bfitb}}_l^{i\top}
\otimes\bigl(\boldintav _l^j - \boldintav
_l^i\bigr)^{\top
}\bigr\} \tilde{\tilde{
\H}}_l \breve{\B}_{l1}(\mathcal{M}) \Biggr]^2
\tilde{w}_{j}^i
\\
& &{} - \sum_{i = 1}^{n}\sum
_{j = 1}^{n} \Biggl[y^j - {
\tilde{a}^i} - \sum_{l = 1}^{g}
\bigl\{{\tilde{\bfitb}}_l^{i\top} \otimes\bigl(\boldintav
_l^j - \boldintav _l^i
\bigr)^{\top}\bigr\} \tilde{\tilde{\H}}_l \tilde{\tilde{
\B}}_{l1} \Biggr]^2 \tilde {w}_{j}^i
\\
& = & 
\sum_{i = 1}^{n}\sum
_{j = 1}^{n} \Biggl[\sum_{l = 1}^{g}
\bigl\{{\tilde {\bfitb}}_l^{i\top} \otimes\bigl(\boldintav
_l^j - \boldintav _l^i
\bigr)^{\top}\bigr\} \tilde{\tilde{\H}}_l \bigl\{ \breve{
\B}_{l1}(\mathcal{M}) - \tilde{\tilde{\B}}_{l1}\bigr\}
\Biggr]^2 \tilde{w}_{j}^i.
\end{eqnarray*}
According to Lemma~4 of Wang and Xia \cite{WanXia08},
there exists some constant
$\kappa> 0$ such that, for any $\mathcal{M} \nsupseteq\mathcal{M}_T$,
\[
\operatorname{RSS}_{\mathcal{M}} - \operatorname{RSS}_{\mathcal{M}_F} \geq\kappa n
h^2
\]
with probability tending to 1. Since $\log(1 + x) \geq\min\{0.5 x,
\log(2)\}$ for any $x > 0$, we have
\begin{eqnarray*}
\log(\operatorname{RSS}_{\mathcal{M}}) - \log(\operatorname{RSS}_{\mathcal{M}_F}) &=& \log
\biggl(1 + \frac{\operatorname{RSS}_{\mathcal{M}} -
\operatorname{RSS}_{\mathcal{M}_F}}{\operatorname{RSS}_{\mathcal{M}_F}} \biggr)
\\
&\geq& \min \biggl\{\log(2), \frac{\operatorname{RSS}_{\mathcal{M}} -
\operatorname{RSS}_{\mathcal{M}_F}}{2\operatorname{RSS}_{\mathcal{M}_F}} \biggr\}.
\end{eqnarray*}
Following an argument similar to the one used in the proof of Lemma~1
of Xia \textit{et al.} \cite{Xiaetal02}, one can show that $n^{-1}\operatorname{RSS}_{\mathcal
{M}_F} \rightarrow_P \sigma^2$ for some $\sigma> 0$. This, together
with $n^{-1} \log(n) = \mathrm{o}(h^2)$, yields that
\[
P \Bigl\{\min_{\mathcal{M} \nsupseteq\mathcal{M}_T} \operatorname{BIC}_{\mathcal{M}} -
\operatorname{BIC}_{\mathcal{M}_F} + \,\mathrm{o}_P
\bigl(h^2\bigr) > 0 \Bigr\} \rightarrow1.
\]
Because $\tilde{\tilde{\B}}_{l1} = \B_{l1}^* + \mathrm{O}_P(n^{-1/2})$ and
$\check{\B}_{l1} = \B_{l1}^* + \mathrm{O}_P(n^{-1/2})$, we obtain
\[
\operatorname{RSS}_{\lambda_n} - \operatorname{RSS}_{\mathcal{M}_F} =
\mathrm{O}_P \biggl(\frac{1}{n} \biggr) =
\mathrm{o}_P\bigl(h^2\bigr).
\]
Thus, we have
\begin{eqnarray*}
P \Bigl(\inf_{\lambda\in\Omega_-} \operatorname{BIC}_{\lambda} -
\operatorname{BIC}_{\lambda_n} > 0 \Bigr) &\geq& P \Bigl(\inf
_{\mathcal{M} \nsupseteq\mathcal{M}_T} \operatorname{BIC}_{\mathcal{M}} -
\operatorname{BIC}_{\mathcal{M}_F} + \operatorname{BIC}_{\mathcal{M}_F} -
\operatorname{BIC}_{\lambda_n} > 0 \Bigr)
\\
&\rightarrow& 1.
\end{eqnarray*}

Next, we consider over-fitted models, that is, $\mathcal{M}_{\lambda}
\supseteq\mathcal{M}_T$ but $\mathcal{M}_{\lambda} \neq\mathcal
{M}_T$. Observe that
\[
\inf_{\lambda\in\Omega_+} \operatorname{RSS}_{\lambda} - \operatorname{RSS}_{\lambda_n} \geq\inf_{\lambda\in\Omega_+} \operatorname{RSS}_{\mathcal{M}_{\lambda}} - \operatorname{RSS}_{\lambda_n} \geq\min
_{\mathcal{M} \supseteq\mathcal{M}_T, \mathcal{M} \neq\mathcal
{M}_T} \operatorname{RSS}_{\mathcal{M}} - \operatorname{RSS}_{\lambda_n}.
\]
For a generic model $\mathcal{M}$, if $\mathcal{M} \supseteq\mathcal
{M}_T$, then one can show that
\[
\operatorname{RSS}_{\mathcal{M}} = \operatorname{RSS}_{\mathcal{M}_F} +
\mathrm{O}_P \biggl(\frac{1}{n} \biggr).
\]
Because $\operatorname{RSS}_{\mathcal{M}_F} - \operatorname{RSS}_{\lambda_n} =
\mathrm{O}_P(n^{-1})$, it follows that
\[
\min_{\mathcal{M} \supseteq\mathcal{M}_T, \mathcal{M} \neq
\mathcal{M}_T} \operatorname{RSS}_{\mathcal{M}} -
\operatorname{RSS}_{\lambda_n} = \mathrm{O}_P \biggl(\frac{1}{n}
\biggr).
\]
Then, with probability tending to 1, we have
\[
\inf_{\lambda\in\Omega_+} \operatorname{RSS}_{\lambda} -
\operatorname{RSS}_{\lambda_n} + \frac{\log(n)}{n} > 0.
\]
As a consequence,
\[
P \Bigl(\inf_{\lambda\in\Omega_+} \operatorname{BIC}_{\lambda} -
\operatorname{BIC}_{\lambda_n} > 0 \Bigr) \geq P \biggl(\inf
_{\lambda\in\Omega
_+} \operatorname{RSS}_{\lambda} - \operatorname{RSS}_{\lambda_n}
+ \frac{\log
(n)}{n} > 0 \biggr) \rightarrow1.
\]
Combining, the proof is complete.
\end{pf*}


%
%
\end{appendix}

\section*{Acknowledgements}

Zhu's research was supported by a
grant from the Research Council of Hong Kong and a grant from Hong Kong
Baptist University,  Hong Kong. Xu's research was supported by the
Natural Science Foundation of Jiangsu Province of China (No.~BK20140617)
and the  Fundamental Research Funds for the Central Universities.


%

\printhistory
\end{document}